\documentclass[12pt]{article}

\usepackage{amssymb}
\usepackage{amsmath}
\usepackage{amsbsy}
\usepackage{amscd}
\usepackage{amsfonts}
\usepackage{amsthm}
\usepackage{mathrsfs}
\usepackage{verbatim}
\usepackage[colorlinks]{hyperref}
\usepackage{fullpage}
\usepackage{mathdots}
\usepackage{graphicx,subfigure}
\usepackage[english]{babel}
\usepackage[utf8]{inputenc}
\usepackage{tikz}

\usetikzlibrary{backgrounds,fit, matrix}
\usetikzlibrary{positioning}
\usetikzlibrary{calc,through,chains}
\usetikzlibrary{arrows,shapes,snakes,automata, petri}
\usepackage{authblk}

\theoremstyle{definition}
\newtheorem{Theorem}{Theorem}[section]
\newtheorem{Corollary}[Theorem]{Corollary}
\newtheorem{Lemma}[Theorem]{Lemma}
\newtheorem{Proposition}[Theorem]{Proposition}
\newtheorem{Example}[Theorem]{Example}

\newtheorem{Remark}[Theorem]{Remark}

\newtheorem{Notation}[Theorem]{Notation}

\newtheorem{Problem}[Theorem]{Open problem}

\newtheorem{Question}[Theorem]{Question}

\newcommand{\C}{{\mathbb C}}
\newcommand{\Z}{{\mathbb Z}}

\newcommand{\mbf}[1]{\mathbf{ #1}}

\newcommand\ind{L}
\newcommand\diag{T}

\newcommand{\BM}{\overline{\mbf{B}}}
\newcommand{\TM}{\overline{\mbf{T}}} 

\newcommand{\Inc}[1]{\mathbf{I}(#1)}
\newcommand{\simtr}{\sim_{\text{tr}}}
\newcommand{\simp}{\sim_{\text{p}}}
\newcommand{\simn}{\sim_{\text{n}}}
\newcommand{\simo}{\sim_{\text{o}}}

\newcommand\gd{\mathcal{D}}
\newcommand\gl{\mathcal{L}}
\newcommand\gr{\mathcal{R}}
\newcommand\gh{\mathcal{H}}
\newcommand\gj{\mathcal{J}}

\begin{document}

\title{Idempotent Varieties of Incidence Monoids and Bipartite Posets}

\author[1]{Mahir Bilen Can}
\author[2,a]{Ana Casimiro}
\author[2,b]{Ant\'onio Malheiro}

\affil[1]{{\small Tulane University, New Orleans LA; \href{mailto:mahirbilencan@gmail.com}{mahirbilencan@gmail.com}}}
\affil[2]{{\small Center for Mathematics and Applications (NOVA Math), and\\ Department of Mathematics, NOVA School of Science and Technology, Portugal}}
\affil[a]{\small\href{mailto:amc@fct.unl.pt}{amc@fct.unl.pt}}
\affil[b]{\small \href{mailto:ajm@fct.unl.pt}{ajm@fct.unl.pt}}


\normalsize

\date{\today}
\maketitle

\begin{abstract}
The algebraic variety defined by the idempotents of an incidence monoid is investigated. 
Its irreducible components are determined. 
The intersection with an antichain submonoid is shown to be the union of these irreducible components.
The antichain monoids of bipartite posets are shown to be orthodox semigroups. 
The Green's relations are explicitly determined, and applications to conjugacy problems are described.
In particular, it is shown that two elements in the antichain monoid are  primarily conjugate in the monoid if and only if 
they belong to the same $\gj$-class and their multiplication by an idempotent of the same $\gj$-class gives conjugate elements in the group.
\medskip

\medskip

\end{abstract}

\section{Introduction}

The incidence monoid of a finite poset is the (complex) linear algebraic monoid whose underlying set consists of ($\C$-valued) 
functions defined on the set of all intervals of the poset, and the multiplication is given by a suitable convolution product. 
In this article we investigate various algebraic subsemigroups in an incidence monoid. 
The purpose of our work is manifold.
First, we show that several important families of semigroups are among the ranks of incidence monoids. 
Secondly, we investigate the structures of the idempotent varieties of such monoids.
Also, we test some notions of conjugacy relations described in~\cite{AKKM} on certain submonoids of the incidence monoids.
Towards achieving these goals, we pay particular attention to the incidence monoids that come from bipartite graphs. The starting point of our analysis is the following theorem. 

\begin{Theorem}\label{T:firstmain}
Let $P$ be a finite poset. 
Let $E(\Inc{P})$ denote the idempotent variety of the incidence monoid of $P$, that is, $E(\Inc{P}):=\{ e\in \Inc{P}:\ e^2 = e \}$.
Then $E(\Inc{P})$ has $2^{|P|}$ connected components. 
Furthermore, each connected component of $E(\Inc{P})$ is an irreducible algebraic subsemigroup of $\Inc{P}$.
\end{Theorem}
Note that our first Theorem~\ref{T:firstmain} does not say that $E(\Inc{P})$ is a subsemigroup, but rather its irreducible components are subsemigroups. In fact, by using the theory of regular algebraic monoids, we will identify a family of posets for which the idempotent varieties are subsemigroups. 
After writing our paper, we learned from Michel Brion that 
a closely related, much more general result about the irreducibility and the smoothness of the components of the idempotent scheme of a not necessarily affine algebraic monoid was already obtained by Brion in~\cite[Theorem 2.14]{Brion2014}. It turns out that our proofs for the relevant parts are different.

The theory of linear algebraic monoids is a fascinating branch of semigroup theory that encompasses the theory of linear algebraic groups. 
An important result of Putcha~\cite{Putcha1984} and Renner~\cite{Renner1985} states that 
the unit group $G(M)$ of an irreducible linear algebraic monoid $M$ with zero is a reductive algebraic group if and only if $M$ is a regular semigroup. 
Here, {\em regular} means that for each $a\in M$ there exists an element $s\in M$ such that $a= asa$. 
When the zero element is missing, the regularity of $M$ is determined by the radical of $G(M)$.
This is given by another theorem of Putcha~\cite{Putcha1984}: an irreducible linear algebraic monoid is regular if and only if the Zariski closure of the radical of $G(M)$ is a completely regular semigroup. Here, {\em completely regular} means that the semigroup is a union of its subgroups. 
This motivates the question of understanding the solvable linear algebraic monoids that are completely regular. 
These monoids are classified in a paper of Renner~\cite{Renner1989}.
We return to the incidence monoid of a finite poset $P$. 
In~\cite{Can2021}, it is shown that there is a set of completely regular submonoids, $\{ \Inc{P,A} \subseteq \Inc{P}:\ \text{$A$ is an antichain in $P$}\}$.
We refer to the elements of this set by {\em antichain monoids (of $P$)}.
Structurally $\Inc{P,A}$ is very similar to $\Inc{P}$.
Indeed, we show in the present article that the idempotent variety of an antichain monoid $\Inc{P,A}$ is a union of certain irreducible components $E(\Inc{P})$ (Proposition~\ref{P:thesameidempotents}).

The antichain monoids of certain posets provide us with important classes of semigroups. 
The second main result of this article is the following. 

\begin{Theorem}\label{T:secondmain}
Let $P$ be a finite poset in which every interval has at most two elements. 
If $A$ is an antichain of $P$, then the corresponding antichain monoid $\Inc{P,A}$ is an {\em orthodox semigroup},
that is, the set of idempotents of $\Inc{P,A}$ is a subsemigroup. 
\end{Theorem}

Let $M$ be an irreducible regular monoid with unit group $G$. 
The {\em cross-section lattice of $M$}, denoted $\Lambda$, is a finite lattice of idempotents of $M$ such that 
$M= \bigsqcup_{e\in \Lambda } GeG$. 
This notion is one of the most important discrete invariants of a regular monoid. 
In general, the computation of this finite poset is difficult. 
Let $P$ be a poset as in Theorem~\ref{T:secondmain}. 
If every maximal element of $P$ covers all minimal elements of $P$, then we call $P$ a {\em complete bipartite poset}.
In the third main result of our paper we analyze the cross-section lattices of the antichain monoids of complete bipartite posets. 

\begin{Theorem}\label{T:thirdmain:intro}
The cross-section lattice of an antichain monoid of a complete bipartite poset is a Boolean lattice. 
\end{Theorem}
 
In an irreducible algebraic monoid $M$ with unit group $G$, a $\gj$-class is given by the two sided orbit $GeG$, where $e$ is an idempotent. 
An important question here is about the structure of the $\gh$-class of $e$. 
What does it look like? 
Our fourth main result answers this question.

\begin{Theorem}\label{T:fourthmain}
Let $e$ be an idempotent of the antichain monoid of a complete bipartite poset $Q$ with respect to the maximal antichain of maximal elements. 
Then the $\gh$-class of $e$ is isomorphic to the group of invertible elements of 
the antichain monoid of a complete bipartite subposet $Q'\subset Q$ with respect to the maximal antichain of maximal elements of $Q'$. 
\end{Theorem}

As a corollary of this result we see that in the antichain monoid $M$ of a complete bipartite poset $Q$ with respect to the maximal antichain $A$ consisting of maximal elements, there are $|A|$ nonisomorphic $\gh$-classes (Corollary~\ref{C:fourthmain}). 
In fact, we determine not only all $\gh$-classes but also every Green's classes for $M$ (Theorem~\ref{T:summaryform}).

As we mentioned earlier, one of our goals in this paper is to initiate a study of conjugacy actions on the incidence monoids.
For a group $G$, there is essentially one type of conjugacy relation; 
$a\sim b$ in $G$ if there exists $g\in G$ such that $a=gbg^{-1}$.
For monoids, there are many different conjugacy notions that agree with the ordinary conjugacy action when restricted to the group of invertible elements. 
To give an example, let us consider a monoid $M$ with unit group $G$. 
The $p$-conjugacy relation is defined by 
\begin{align*}
a\simp b \Leftrightarrow \exists z,\, w\in M:\, a=zw ,\, b=wz.
\end{align*}
It is easy to see that on $G$, we have $\sim\ \equiv \ \simp$. 
In the last main result of our paper, 
we determine the $\simp$ conjugacy classes in the antichain monoid of a complete bipartite poset.

\begin{Theorem}\label{T:lastannounced}
Let $M$ denote the antichain monoid (associated with a maximal antichain) of a complete bipartite poset. 
Let $X$ and $Y$ be two elements from $M$. 
Then $X \simp Y$ if and only if 
both $X$ and $Y$ belong to the same $\gj$-class and their (right) multiplication by an idempotent of the same $\gj$-class gives conjugate elements in the group.
\end{Theorem}

Let us finally mention that this theorem shows that the $p$-conjugacy problem in $M$ is reduced to the ordinary conjugacy problem in the $\gh$-classes. We solve this problem explicitly in our Corollary~\ref{C:moreprecisely}.
\medskip

We are now ready to describe the contents of the individual sections.
In the next preliminaries section, we setup our notation, and review the fundamentals of linear algebraic monoids. 
Also in this section, we introduce the incidence monoids and their antichain submonoids. 
The purpose of Section~\ref{S:Idempotents} is to analyze the idempotent variety of an incidence monoid. 
This is where we prove our first main result Theorem~\ref{T:firstmain}.
In addition, we show that the components of the idempotent variety of an antichain submonoid are among the irreducible components of the idempotent variety of the ambient incidence monoid (Proposition~\ref{P:thesameidempotents}). 
At the beginning of Section~\ref{S:Bipartite}, we prove our Theorem~\ref{T:secondmain}, which states roughly that the idempotents of the incidence monoid of a bipartite poset is a subsemigroup. In Theorem~\ref{T:componentmultiplication}, we analyze the product of two irreducible components in this semigroup of idempotents. 
In Section~\ref{S:Cross}, we analyze the lattice of regular $\gj$-classes of the antichain monoid of a complete bipartite poset. 
More precisely, we prove our Theorem~\ref{T:thirdmain:intro}, which states that the lattice of regular $\gj$-classes here form a Boolean lattice.
Another goal of Section~\ref{S:Cross} is to describe the $\gh$-classes of the idempotents. 
It turns out that the $\gh$-classes in our particular antichain monoids look like the unit groups of appropriate antichain monoids.
This is where we prove our Theorem~\ref{T:fourthmain}. 
In Section~\ref{S:Green}, we describe explicitly all Green equivalence classes of an arbitrary element of a maximal antichain monoid of a complete bipartite poset (Theorem~\ref{T:summaryform}).
This result turns out to be very instrumental for our study of the $\simp$-conjugacy on the antichain monoid of a complete bipartite poset.
We prove our Theorem~\ref{T:lastannounced} in Section~\ref{S:Conjugacy}.
Finally, we close our paper by a brief discussion of another equivalence relation that is closely related to $\simp$. 
We pose several open problems for future research.

\section{Preliminaries}\label{S:Preliminaries}

Although our results hold true for any algebraically closed field of characteristic zero, for simplicity, we work with algebraic semigroups defined over the field of complex numbers.
We denote the set of positive integers by $\Z_+$. 
Let $\{k,n\}\subset \Z_+$. 
We fix the following notation for the rest of our text:
\begin{align*}
\begin{array}{ccl}
\mbf{Mat}_{k,n} &:& \text{the set of $k\times n$ matrices};\\
\mbf{Mat}_n &:& \text{the monoid of $n\times n$ matrices};\\
\mbf{GL}_n &:& \text{the unit group of $\mbf{Mat}_n$};\\
\mbf{B}_n &:& \text{the subgroup of upper triangular matrices in $\mbf{GL}_n$};\\
\mbf{U}_n &:& \text{the subgroup of unipotent upper triangular matrices in $\mbf{B}_n$};\\
\mbf{T}_n &:& \text{the diagonal torus, that is, the subgroup of diagonal matrices in $\mbf{B}_n$};\\
\mbf{1}_n &:& \text{the $n\times n$ identity matrix};\\
\text{$[n]$} &:& \text{the set $\{1,\dots, n\}$};\\
\text{diag}(a_1,\dots, a_n) &:& \text{$n\times n$ diagonal matrix with entries $a_1,\dots, a_n$}.\\
\end{array}
\end{align*}
For $\ind \subseteq [n]$, we denote by $\mbf{1}_\ind$ the diagonal 
idempotent
\begin{align}\label{A:1_J}
\mbf{1}_\ind:= \text{diag}(a_1,\dots, a_n),
\ \text{ where } 
a_j =
\begin{cases}
1 & \text{ if } j\in \ind\\
0 & \text{ otherwise}.
\end{cases}
\end{align}

\medskip

Let $G$ be an algebraic group. 
The {\em radical} of $G$, denoted $R(G)$, is the maximal connected normal solvable subgroup of $G$. 
If $R(G)$ is trivial, then $G$ is said to be {\em semisimple}. 
The {\em unipotent radical} of $G$, denoted $R_u(G)$, is the maximal connected normal unipotent subgroup of $G$. 
If $R_u(G)$ is trivial, then $G$ is said to be {\em reductive}.

\medskip

An {\em algebraic semigroup} is an algebraic variety $S$ with an associative binary operation 
$\mu: S\times S\to S$ such that $\mu$ is a morphism of varieties. 
Following the useful conventions of algebraic group theory, in this text, 
we do not assume that varieties are irreducible.
Of course, an algebraic set has only finitely many irreducible components.  
An algebraic monoid is an algebraic semigroup with an identity element. 
For an introduction to the theory of not necessarily linear algebraic semigroups, we refer the reader
to~\cite{Brion:FieldsLectures,Renner}.

\medskip

We now introduce some poset theory terminology. 
Let $(P,\leq)$ be a poset. 
We will usually omit the order relation ``$\leq$'' from our notation. 
A {\em chain} is a totally ordered subposet of $P$. 
A {\em maximal chain} in $P$ is a chain $C$ in $P$ such that $C$ is not a chain in any other chain in $P$.  Let $P$ be a finite poset. 
If $C$ is a chain, then its {\em length} is the number $|C|-1$.
If every maximal chain in $P$ has the same length, then $P$ is said to be {\em ranked} (or, {\em graded}).
In this case, the length of any maximal chain in $P$ is called the {\em rank of $P$}.
In this paper we are concerned with the incidence monoids of finite posets only.

\subsection{Semigroups, algebraic semigroups.}\label{SS:Preliminaries:Semigroups}

The group $\mbf{B}_n$ is a Borel subgroup of $\mbf{GL}_n$. 
We call the Zariski closure of $\mbf{B}_n$ in $\mbf{Mat}_n$ the {\em (standard) Borel monoid}.
The main purpose of this subsection is to review some general results of Putcha on the closed subsemigroups of the Borel monoid. 
We begin with some general semigroup theory notions. We shall use the following results without further notice. For more information on semigroups see~\cite{Howie}.
\medskip

Throughout this section $S$ will denote a semigroup unless otherwise specified. 
Let $S^1$ denote the monoid obtained by adjoining an identity to $S$ if necessary (in which case $S^1 := S\cup \{1\}$). 
The {\em Green's relations on $S$} are defined as follows.
Let $a,b\in S$. Then 
\begin{enumerate}
\item $a\gr b$ if $aS^1 = bS^1$,
\item $a\gl b$ if $S^1a = S^1b$,
\item $a\gj b$ if $S^1aS^1 = S^1bS^1$,
\item $a\gh b$ if $a\gr b$ and $a\gl b$,
\item $a\gd b$ if $a\gr c$ and $c\gl b$ for some $c\in S$.
\end{enumerate}
The equivalence class of an element $a\in S$ with respect to any of these relations will be indicated by putting $a$ in the subscript. 
For example, $H_a$ stands for the $\gh$-class of $a$. 
\medskip

We denote the set of idempotents of $S$ by $E(S)$.
The natural partial order on $E(S)$ is given by 
\begin{align}\label{A:naturalorderonE(S)}
e\leq f \iff e=fe=ef
\end{align}
for $e,f\in E(S)$. 
Let $I$ be a nonempty subset of $S$.
Then $I$ is called a {\em right ideal of $S$} (resp. {\em left ideal}, resp. {\em ideal}) if $IS \subseteq I$ (resp. $SI\subseteq I$, resp. $S^1 I S^1 \subseteq I$). The semigroup $S$ is called {\em simple} if it contains no proper ideals, and if $S$ contains a zero element $0$, then it is said {\em $0$-simple} if it is not a null semigroup and its only proper ideal is 
$\{0\}$.



For $a\in S$, we define the set 
$V(a) := \{ s\in S :\ asa=a\text{ and } s=sas\}$.
We call an element $s$ of $V(a)$ an {\em inverse} of $a$.
Many useful definitions of semigroup theory can be stated via the set  $V(a)$.
For example, a {\em regular semigroup} is a semigroup $S$ such that  $V(a) \neq \emptyset$ for every $a\in S$. 
\medskip

Recall that $S$ is called {\em completely regular} if it is a union of subgroups. 
There are several useful descriptions of a completely regular semigroup. 
Here is one of them: 
$S$ is completely regular if and only if for every $a\in S$ there is an inverse $a^{-1} \in V(a)$ such that 
the following identities hold:
\begin{align}\label{A:inverseincrs}
(a^{-1})^{-1}=a \ \text{ and }\  a^{-1}a=aa^{-1}.
\end{align}
Note that in a completely regular semigroup an inverse $a^{-1} \in V(a)$ satisfying the conditions in (\ref{A:inverseincrs}) is unique.

The following facts will be useful for our purposes.
\begin{enumerate}
\item Any regular subsemigroup of $\BM_n$ is completely regular~\cite[Remark 3.21]{Putcha}.
\item Let $M$ be an irreducible linear algebraic monoid with unit group $G$. Let $e$ be a minimal idempotent of $M$. 
Then $M$ is a regular semigroup if and only if $G_e:= \{g\in G:\ ge =eg =e \}$ is a reductive group~\cite[Theorem 7.4]{Putcha}.
\end{enumerate}
{\it Warning:} In Putcha's monograph~\cite{Putcha}, a connected monoid means an irreducible monoid.
\medskip

A fundamental theorem due to Putcha~\cite[Theorem 3.18]{Putcha} for linear algebraic semigroups, Brion and Renner for arbitrary algebraic semigroups~\cite{BrionRenner}, states that every algebraic semigroup is {\em strongly $\pi$-regular}, also known as {\em epigroup}, that is, for any $x\in S$, there exists a positive integer $n\in \Z_+$ such that 
$x^n$ lies in a group. 
Let $S$ be an epigroup. 
Let $a$ and $e$ be two elements of $S$ such that $e\in E(S)$. 
Proofs of the following statements can be found in~\cite[Theorem 2.65]{Renner}:
\begin{enumerate}
\item Let $J_e \cap e S e = H_e$;
\item $\gj=\gd$ on $S$; 
\item If $a\gj a^2$, then $H_a$ is a group. 
In particular, $H_e$ is an algebraic group. 
\end{enumerate}
A $\gj$-class $J$ of $S$ is said to be {\em regular} if some (hence every) element of $J$ is regular.  
This is equivalent to the requirement that  $E(J)\neq \emptyset$.
We will denote by $\mathscr{U}(S)$ the partially ordered set of all regular $\gj$-classes of $S$. 
Here, the partial order on $\mathscr{U}(S)$ is defined as follows: for any regular elements $a$ and $b$ in $S$
\begin{align}\label{A:partialorderonJ}
J_b \leq J_a \iff \text{ $xay=b$ for some $x,y\in S^{1}$}.
\end{align}
In~\cite[Theorem 3.28]{Putcha} Putcha shows that if $S$ is an algebraic semigroup, then $\mathscr{U}(S)$ is a finite poset.

\medskip


A semigroup $S$ is said {\em completely simple} if it is a simple epigroup~\cite[Theorem 3.3.2]{Howie}, and it is said {\em completely 0-simple} if it is a $0$-simple epigroup \cite[Theorem 3.2.11]{Howie}. 
It is pointed out in~\cite[Remark 2.68]{Renner} that these are not the standard definitions of completely simple and completely $0$-simple semigroups.  
However, by a theorem of Munn, they are equivalent to the standard ones.
A semigroup $S$ is called {\em right simple} if the $\gr =S\times S$.
A semigroup which is right simple and left cancellative is called a {\em right group} \cite[pag. 61]{Howie}.
A famous result of Clifford from~\cite{Clifford} states that $S$ is completely regular if and only if $S$ is a semilattice of completely simple semigroups (c.f. \cite[\S 4.1]{Howie}). In this case, $S/\gj$ is a semilattice and each $\gj$-class is a completely simple semigroup.


\medskip 

In this paragraph, $M$ denotes an irreducible (hence connected) linear algebraic monoid with unit group $G$.
Let $\{a,b\}\subset M$ and $\{e,f\}\subset E(M)$. 
Since $G$ is dense in $M$, we know the following facts from~\cite[Proposition 6.1]{Putcha}:
\begin{enumerate}
\item $a\gr b$ if and only if $aG = bG$;
\item $a\gl b$ if and only if $Ga = Gb$;
\item $a\gj b$ if and only if $GaG = GbG$.
\end{enumerate}
We know the following facts from~\cite[Proposition 6.8]{Putcha}:
\begin{enumerate}
\item $e\gj f$ if and only if $x^{-1}ex =f$ for some $x\in G$;
\item $e\gr f$ if and only if there exists $x\in G$ such that $ex=x^{-1}ex=f$;
\item $e\gl f$ if and only if there exists $x\in G$ such that $xe=xex^{-1} = f$.
\end{enumerate}
Finally, we state the structure theorem~\cite[Corollary 6.10]{Putcha} for the idempotent variety of $M$.
For any maximal torus $T$ in $G$, the idempotent variety of $M$ is given by the union 
\begin{align}\label{A:Minconjugacy}
E(M) = \bigcup_{x\in G} x^{-1} E(\overline{T})x.
\end{align}

\medskip

We are now ready to review Putcha's ``cross-section lattice''.
Let $M$ be an irreducible linear algebraic monoid with unit group $G$. 
A subset $\Lambda$ of $E(M)$ is called a {\em weak cross-section lattice} if 
\begin{itemize}
\item $|\Lambda \cap J| =1 $ for all $J\in \mathscr{U}(M)$, and 
\item If $e,f\in \Lambda$, then $J_e \leq J_f$ implies $e\leq f$. 
\end{itemize}
If in addition $\Lambda \subseteq E(\overline{T})$ for some maximal torus $T$ of $G$, then $\Lambda$ is called a {\em cross-section lattice}.
Since the regular $\gj$-classes of an irreducible linear algebraic semigroup form a finite lattice (\cite[Theorem 5.10]{Putcha}), the cross-section lattice of an irreducible linear algebraic monoid is indeed a finite lattice. In particular, every cross-section lattice has a unique minimal and a unique maximal elements.

Recall that the $\gj$-class of an element $x\in M$ 
is given by $GxG$ (\cite[Proposition 6.1]{Putcha}). Thus, it follows from the definition of $\Lambda$ that the natural partial order inherited from $E(M)$ on $\Lambda$ agrees with the following order: $e\leq f \iff GeG \subseteq \overline{GfG}$ for every $\{e,f\}\subseteq \Lambda$.

\medskip

We close this section by summarizing some important properties of the regular linear algebraic monoids.
If $M$ is an irreducible regular monoid, then the following statements hold: 
\begin{enumerate}
\item $M$ has a cross-section lattice~\cite[Corollary 9.4]{Putcha};
\item Any two cross-section lattices of $M$ are conjugate~\cite[Corollary 9.7]{Putcha}.
\end{enumerate}

\subsection{Incidence monoids.}

Incidence monoids come from certain associative algebras. 
More information on such semigroups can be found in Okninski's article,~\cite{Okninski2014}.
\medskip

Let $(P,\leq)$ be a poset.
The set of all intervals of $P$ is denoted by $\text{int}(P)$. 
Let $f$ and $g$ be two $\C$-valued functions on $\text{int}(P)$. 
The {\em convolution product of $f$ and $g$} is the product defined by 
\begin{align*}
(f* g) ( [s,u])= \sum_{s\leq t \leq u } f([s,t]) g([t,u]) \qquad ([s,u]\in \text{int}(P)).
\end{align*}
The {\em incidence monoid of $P$}, denoted $\Inc{P}$, is the monoid of $\C$-valued functions on $\text{int}(P)$,  
where the multiplication is given by the convolution product.
There is a vector space structure on $\Inc{P}$, where the sum is given by the point-wise addition of the functions. 
It follows that, if $P$ is a finite poset, then $\Inc{P}$ is an affine space.
In particular, in this case, $\Inc{P}$ is a linear algebraic monoid.
Hereafter, we will work with finite posets only.
\medskip

We proceed with descriptions of certain linear representations of $\Inc{P}$.
A {\em linear extension of $P$} is a bijection $\mathbf{le}: P\to [n]$ such that, for every $\{s,u\}\subset P$,
the following implication holds: 
\begin{align*}
s\leq u \implies \mathbf{le}(s) \leq \mathbf{le}(u). 
\end{align*}
Here, $\leq$ is the natural order on $[n]$. 
From now on, when we list the elements of $P$ as in $P=\{x_1,\dots,x_n\}$, 
we assume implicitly that the assignment $x_i \mapsto i$, $i\in [n]$ is a fixed linear extension of $P$. 
For such a presentation of $P$, we get a $\C$-algebra representation 
\begin{align*}
\Psi=\Psi_{\mbf{le}}: \Inc{P} &\longrightarrow \mbf{Mat}_n\\
f &\longmapsto (f([x_i,x_j]))_{i,j=1}^n,
\end{align*}
where we set $f([x_i,x_j]):=0$ for every non-relation $x_i \nleqslant x_j$ in $P$.
Clearly, the image of $\Psi$ is a closed submonoid of $\BM_n$.
Indeed, $\Psi$ is a faithful representation of $\Inc{P}$. 
Also, since any two matrix representations that are obtained from two (different) linear extensions of $P$ are isomorphic, 
for our purposes, the choice of a linear extension is unimportant.  
Hereafter, we will identify the incidence monoid $\Inc{P}$ with its image $\Psi(\Inc{P})$ in $\BM_n$.

\begin{Example}\label{E:Borelmonoid}
Let $P$ be a chain with $n$ elements.
In other words, let $P=\{x_1,\dots,x_n\}$, where $x_1< \cdots < x_n$.
The image of the corresponding matrix representation of $\Inc{P}$ is the Borel monoid, $\BM_n$. 
\end{Example}

The unit group of an incidence monoid $\Inc{P}$ will be denoted by $\mbf{G}(P)$.
We note that since $\Inc{P}$ has the structure of a linear space, it is irreducible as a variety. 
In particular, the unit group $\mbf{G}(P)$ is an irreducible algebraic group.

\subsection{Antichain monoids.}

We are now ready to properly introduce the antichain monoids. 
Let $P:=\{x_1,\dots, x_n\}$ be a poset. 
Let $A$ be a nonempty antichain in $P$. 
The {\em antichain monoid}, denoted $\Inc{P,A}$, is the linear algebraic monoid defined by the Zariski closure $\overline{\mbf{T}_n(A)\ltimes \mbf{U}_n(P)} \subseteq \mbf{Mat}_n$, where 
$\mbf{T}_n(A)$ is the diagonal torus, 
\begin{align}\label{A:ourdiagonaltorus}
\mbf{T}_n(A) := \left\{ \text{diag}(t_1,\dots, t_n) \in \mbf{Mat}_n :\ \text{$t_i =1$ if $x_i\notin A$ for $i\in [n]$}\right\}, 
\end{align} 
and $\mbf{U}_n(P)$ is the group of upper triangular unipotent matrices $(a_{ij})_{i,j=1}^n \in \mbf{Mat}_n$ such that $a_{ij}=0$ for every $i,j\in \{1,\dots, n\}$ whenever $x_i \nleq x_j$.
Note that $\mbf{U}_n(P)$ is the unipotent radical of the unit group of $\Inc{P}$.
We denote the unit group of $\Inc{P,A}$ by $\mbf{G}(P,A)$.
Since we have $\mbf{G}(P,A) = \mbf{T}_n(A)\ltimes \mbf{U}_n(P)$, the antichain monoid $\Inc{P,A}$ is a unit-dense monoid. In~\cite[Proposition 2.4]{Can2021} it is shown that $\Inc{P,A}$ is a completely regular semigroup.
\medskip

Let us call a poset {\em connected} if its Hasse diagram is connected.
A connected component of a finite poset $P$ is a subposet $P'\subseteq P$ such that the Hasse diagram of $P'$ is a connected component of the Hasse diagram of $P$. 
The following proposition will be useful in the sequel. 
\begin{Proposition}\label{P:connectedcomponents}
Let $A$ be an antichain of $P$. 
Let $P_1,\dots, P_s$ denote the connected components of $P$. 
If $A_i$ ($i\in [s]$) denotes the antichain defined by $A_i:=A\cap P_i$, then we have the following decomposition of algebraic monoids:
\begin{align*}
\Inc{P,A} = \mbf{I}(P_1,A_1) \times \cdots \times \mbf{I}(P_s,A_s).
\end{align*}
\end{Proposition}
\begin{proof}
For $i\in [n]$, let $n_i$ denote the cardinality of $P_i$. 
The proofs of the following decompositions follow from definitions, 
\begin{align*}
\mbf{U}_n(P) = \mbf{U}_{n_1}(P) \times \cdots \times \mbf{U}_{n_s}(P)
\qquad \text{and} \qquad \mbf{T}_n(A) = \mbf{T}_{n_1}(A_1) \times \cdots \times \mbf{T}_{n_s}(A_s).
\end{align*} 
It is now easy to see that 
$\mbf{G}(P,A) =  \prod_{i=1}^s (\mbf{T}_{n_i}(A_i)\ltimes \mbf{U}_{n_i}(P_i) )$.
Since the Zariski closure of the right hand side of this isomorphism is the product of the Zariski closures of the factors, the proof of our assertion follows.
\end{proof}

\section{The Idempotent Varieties}\label{S:Idempotents}

In this section we investigate the structure of the idempotent variety of an incidence monoid. 
As we mentioned before, the idempotent variety of a not necessarily linear algebraic semigroup was investigated by Brion in~\cite{Brion2014}.
Some of our general statements in this section can be derived from Brion's work. 
For completeness, we will provide the proofs of all of our statements. 
Let us begin with a basic example to motivate our discussion. 

\begin{Example}
Let $Q$ denote the chain $x_1< x_2$. 
The incidence monoid of $Q$ is the Borel monoid $\BM_2$.
The idempotent variety of $Q$, that is $E(Q)$, has four connected components. 
They are given by 
\begin{enumerate}
\item $C_1 : = \left\{  \begin{bmatrix} 1 & a \\ 0 & 0 \end{bmatrix} :\ a\in \C   \right\}$,
\item $C_2 : = \left\{  \begin{bmatrix} 0 & b \\ 0 & 1 \end{bmatrix} :\ b\in \C   \right\}$, 
\item $C_3 :=  \left\{  \begin{bmatrix} 1 & 0 \\ 0 & 1 \end{bmatrix} \right\}$, 
\item $C_4 :=  \left\{  \begin{bmatrix} 0 & 0 \\ 0 & 0 \end{bmatrix} \right\}$.
\end{enumerate}
Clearly, $Q$ has two antichains, $A:=\{x_2\}$ and $A':=\{x_1\}$. 
It is not difficult to verify that the union $C_1\cup C_3$ is the idempotent variety of $\Inc{Q,A}$,
and the union $C_2\cup C_3$ is the idempotent variety of $\mbf{I}(Q,A')$.
\end{Example}

Let $P=\{x_1,\dots, x_n\}$ be a poset. 
The diagonal monoid $\TM_n$ (the Zariski closure of $\mbf{T}_n$) is a maximal abelian algebraic submonoid of $\Inc{P}$.
The following natural map of upper triangular matrices is an algebraic monoid homomorphism between $\BM_n$ and $\TM_n$,
\begin{align}\label{A:naturalmapofupper}
(a_{ij})_{i,j=1,\dots,n} \mapsto \text{diag}(a_{11},\dots, a_{nn})\qquad ((a_{ij})_{i,j=1}^n\in \BM_n).
\end{align}
The restriction of (\ref{A:naturalmapofupper}) to its closed submonoid $\Inc{P}$ is still a surjective algebraic monoid homomorphism. 
Thus, by~\cite[Corollary 1]{Brion:FieldsLectures}, we get a morphism of idempotent varieties,
\begin{align}
p: E(\Inc{P}) &\longrightarrow E(\TM_n)\notag \\
(e_{ij})_{i,j=1,\dots, n} &\longmapsto \text{diag}(e_{11},\dots, e_{nn}).\label{A:EtoE}
\end{align}
Our final remark before proving the main result  of this section is the following:
Since $\TM_n$ is abelian, the variety $E(\TM_n)$ is a finite set of points,~\cite[Proposition 4 (iii)]{Brion:FieldsLectures}.
In fact, $E(\TM_n)$ is given by 
\begin{align*}
E(\TM_n) = \{ \mbf{1}_\ind :\ \ind\subseteq [n]\},
\end{align*}
where $\mbf{1}_\ind$ is the idempotent defined in (\ref{A:1_J}).

We are now ready to prove the first announced theorem of our paper. 
Let us recall its statement for convenience. 
\begin{Theorem}\label{T:variety_idempotents}
Let $P=\{x_1,\dots, x_n\}$ be a poset. 
The idempotent variety $E(\Inc{P})$ has $2^n$ connected components.
Each connected component is of the form
\begin{align}\label{A:twoequations}
p^{-1}(\mbf{1}_\ind)=\mbf{G}(P)\cdot \mbf{1}_\ind =J_{\mbf{1}_\ind}\cap E(\Inc{P}), 
\end{align}
for some $\ind\subseteq [n]$. 
Furthermore, each connected component is an irreducible subsemigroup of $\Inc{P}$.
\end{Theorem}

In (\ref{A:twoequations}), the action of $\mbf{G}(P)$ is the conjugation action; 
$J_{\mbf{1}_\ind}$ stands for the $\gj$-class of ${\mbf{1}_\ind}$.

\begin{proof}
Since $\TM_n$ is a submonoid of $\Inc{P}$, the morphism $p$ in (\ref{A:EtoE}) is surjective. 
The image of $p$ has $2^n$ elements. 
Let $\ind\subseteq [n]$. 
We will analyze the subvariety $p^{-1}(\mbf{1}_\ind)$ of $E(\Inc{P})$.

For $\{k,l\}\subseteq [n]$, let $x_{kl}$ denote the $(k,l)$-th coordinate function on $\mathbf{Mat}_n$. 
Recall that $\Inc{P}$ is the affine subspace of $\mathbf{Mat}_n$ defined by the equations $x_{kl}=0$ corresponding to the non-relations in $P$. 
In this notation, the closed subset $p^{-1}(\mbf{1}_\ind) \subseteq \Inc{P}$ is defined by the conditions,  
\begin{align*}
x_{jj} = 0 \text{ for $j\in [n] \setminus \ind$},\qquad 
x_{ii} = 1 \text{ for $i\in \ind$},
\end{align*}
and 
\begin{align}\label{A:thesecondsetofeqns}
\sum_{l=1}^n x_{s l} x_{l t} &= x_{st} \text{ for $\{s,t\}\subseteq [n]$}.
\end{align}
Since the equations in (\ref{A:thesecondsetofeqns}) do not have constant terms, by sending the off-diagonal matrix coordinate variables to 0, we see the idempotent $\mbf{1}_\ind$ is contained in every (path-)connected component of the preimage, $p^{-1}(\mbf{1}_\ind)$. 
In other words, $p^{-1}(\mbf{1}_\ind)$ is a Zariski connected topological space. 
Therefore, $E(\Inc{P})$ has $2^n$ connected components. 
This finishes the proof of our first assertion.

We proceed to show that the connected component $p^{-1}(\mbf{1}_\ind)$ is irreducible. 
Recall that $\Inc{P}$ is irreducible and that $E(\TM_n)=\{\mbf{1}_\ind:\ \ind\subseteq [n]\}$. 
Then $\mbf{G}(P)$ is irreducible, and hence the orbit $\mbf{G}(P)\cdot \mbf{1}_L$ is an irreducible subvariety of $E(\Inc{P})$. 
From the decomposition (\ref{A:Minconjugacy}) it follows that 
\begin{align*}
E(\Inc{P}) = \bigcup_{\mbf{1}_L \in E(\TM_n)} \mbf{G}(P)\cdot \mbf{1}_L,
\end{align*}
where $\mbf{G}(P)\cdot \mbf{1}_L = \{ x^{-1} \mbf{1}_L x :\ x\in \mbf{G}(P)\}$. 
Since $\mbf{G}(P) \cdot \mbf{1}_L \subseteq p^{-1}(\mbf{1}_L)$, the above decomposition is disjoint.  
It follows immediately that $p^{-1}(\mbf{1}_L) = \mbf{G}(P)\cdot \mbf{1}_L$, and hence is irreducible.  
Also notice that any idempotents are in the $\mathcal{J}$-class if and only if both belong to the same conjugacy orbit by~\cite[Proposition 6.8]{Putcha}.

Finally, if $X$ and $Y$ are two idempotents from $p^{-1}(\mbf{1}_\ind)$ for some $\ind\subseteq [n]$, 
then the diagonal of $XY$ is equal to the diagonal of $X$. In other words,
$p(XY)=\mbf{1}_\ind$. It follows that $p^{-1}(\mbf{1}_\ind)$ is closed under multiplication.
Hence, it is a subsemigroup.  
This finishes the proof of our theorem. 
\end{proof}

\begin{Remark}
In a personal communication, it is pointed to us by Michel Brion that the second part of the conclusion of our theorem holds more generally.
Let $M$ be an irreducible algebraic monoid with unit group $G$. 
It is proved in~\cite[Theorem 2.1]{Brion2014} that the scheme of idempotents $E(M)$ is smooth. 
Furthermore, each connected component of $E(M)$ is a conjugacy class of $G$; it meets the closure of a given maximal torus of $G$. 
Consequently, each connected component of $E(M)$ is irreducible. 
\end{Remark}

\begin{Example}
Let $P$ denote the chain $x_1 < x_2 < x_3$. 
Then we have $\Inc{P}= \BM_3$. 
Let us consider the irreducible component $p^{-1}(\mbf{1}_\ind)$ of $E(\BM_3)$, 
where $\ind$ is the subset $\ind:=\{1,3\}\subset [3]$. 
It is easy to check that 
\begin{align*}
p^{-1}(\mbf{1}_\ind) = 
\left\{ 
\begin{bmatrix}
1 & a & b \\
0 & 0 & c \\
0 & 0 & 1
\end{bmatrix}
: b= -ac
\right\}.
\end{align*}
Clearly, this is isomorphic as a variety to the irreducible hypersurface defined by the vanishing of the polynomial $y+xz$ in the affine $3$-space, $\mathbb{A}^3(\C)$. 
\end{Example}

The following observation will be useful in the sequel. 
\begin{Proposition}\label{P:thesameidempotents}
Let $\Inc{P,A}$ be an antichain submonoid of $\Inc{P}$ for some antichain $A\subseteq P$.
Then every irreducible component of $E(\Inc{P,A})$ is an irreducible component of $E(\Inc{P})$.
In particular, if $\mbf{1}_\ind$ is a diagonal idempotent of $\Inc{P,A}$ for some $\ind\subseteq [n]$, then 
we have 
\begin{align}\label{A:twoequalitiesforG(P,A)}
p^{-1}(\mbf{1}_\ind)=\mbf{G}(P,A)\cdot \mbf{1}_\ind =J_{\mbf{1}_\ind}(A)\cap E(\Inc{P,A}),
\end{align}
where $J_{\mbf{1}_\ind}(A)$ is the $\gj$-class of $\mbf{1}_\ind$ in $\Inc{P,A}$.
\end{Proposition}

\begin{proof} 
Since $E(\Inc{P,A}) = E(\Inc{P})\cap \Inc{P,A}$, by Theorem~\ref{T:variety_idempotents}, 
it suffices to show that $p^{-1}(\mbf{1}_\ind)\subseteq E(\Inc{P,A})$ for every $\ind\subseteq [n]$. 
Of course, if for some $\ind\subseteq [n]$ the idempotent $\mbf{1}_\ind$ is not an element of $E(\Inc{P,A})$, then there is nothing to prove. So, we proceed with the assumption that $\mbf{1}_\ind \in E(\Inc{P,A})$ for some subset $\ind\subseteq [n]$.
Let $x$ be an idempotent from $p^{-1}(\mbf{1}_\ind)$.
Then $x$ is of the form $x= \mbf{1}_\ind + N$, where $N$ is a strictly upper triangular $n\times n$ matrix.
We now define a one-parameter subgroup of $\mbf{G}(P,A)$:
\begin{align}
\lambda : \C^* &\longrightarrow \mbf{G}(P,A) \notag \\
t &\longmapsto \mbf{1}_\ind + t(\mbf{1}_n-\mbf{1}_L) + N \label{A:1PSG}
\end{align}
Clearly, we have $\lim_{t\to 0} \lambda(t) = \mbf{1}_\ind + N$.
Hence, the idempotent $\mbf{1}_\ind + N$ is contained in the closure of the image of $\lambda$.
At the same time, the following inclusion holds: 
\begin{align*}
\lambda (\C^*)\subseteq \overline{\mbf{G}(P,A)}=\Inc{P,A}.
\end{align*}
Therefore, the idempotent $\mbf{1}_\ind+N$ is contained in $\Inc{P,A}$. 
This finishes the proof of our first assertion.

To prove our second assertion, we recall the decomposition (\ref{A:Minconjugacy}):
the idempotents of an irreducible monoid $M$ with unit group $G$ 
are given by the $G$-conjugates of the idempotents of the closure of a maximal torus of $G$.
A maximal torus of $\mbf{G}(P,A)$ is given by the diagonal torus $\mbf{T}_n(A)$ as defined in (\ref{A:ourdiagonaltorus}). 
Since every element of $\mbf{G}(P,A)$ is upper triangular, 
two-idempotents $\ind$ and $\ind'$ of $\mbf{T}_n(A)$ are $\mbf{G}(P,A)$-conjugate if and only if they are equal. 
It follows that $\mbf{G}(P,A)\cdot \mbf{1}_\ind$ does not contain any other idempotent of $\mbf{T}_n(A)$. 
Since $\mbf{G}(P,A)$ is connected, so is its orbit $\mbf{G}(P,A)\cdot \mbf{1}_\ind$.
This means that $\mbf{G}(P,A)\cdot \mbf{1}_\ind$ is a connected component of $E(\Inc{P,A})$. 
Since $ \mbf{G}(P,A)\cdot \mbf{1}_\ind \subseteq \mbf{G}(P)\cdot \mbf{1}_\ind$, and $\mbf{G}(P)\cdot \mbf{1}_\ind$ is a connected component of 
$E(\Inc{P,A})$ also, we see that $ \mbf{G}(P,A)\cdot \mbf{1}_\ind = \mbf{G}(P)\cdot \mbf{1}_\ind$.
Hence, the proof of the first equality in (\ref{A:twoequalitiesforG(P,A)}) follows.
The proof of the second equality is similar, so, we skip it. 
This finishes the proof of our proposition.
\end{proof}

\begin{Remark}\label{R:twoitems}
Let $A$ and $A'$ be two antichains in $P$. 
By~\cite[Corollary 2.5]{Can2021} we know that 
\begin{enumerate}
\item $A'\subseteq A$ if and only if $\Inc{P,A'}\subseteq \Inc{P,A}$, and 
\item $\Inc{P,A\cap A'}=\Inc{P,A}\cap \Inc{P,A'}$.
\end{enumerate}
These items combined with Proposition~\ref{P:thesameidempotents} show that 
\begin{enumerate}
\item $A'\subseteq A \Rightarrow E(\Inc{P,A'})\subseteq E(\Inc{P,A})$, and 
\item $E(\Inc{P,A\cap A'})=E(\Inc{P,A})\cap E(\Inc{P,A'})$.
\end{enumerate}
\end{Remark}

We close this section by another useful observation.
\begin{Lemma}\label{L:smooth}
Let $A$ be an antichain of $P:=\{x_1,\dots, x_n\}$. 
Then $\Inc{P,A}$ is an affine subspace of $\BM_n$. 
\end{Lemma}
\begin{proof}
Let $\mbf{1}_A$ denote the minimal diagonal idempotent of $\Inc{P,A}$.
In other words, we have $\mbf{1}_A:=\text{diag}(a_1,\dots, a_n)$, where
\begin{align*}
a_i :=
\begin{cases}
0 & \text{ if $x_i \in A$}\\
1 & \text{ if $x_i \notin A$}.
\end{cases}
\end{align*}
We observe that $\overline{\mbf{T}_n(A)}- \mbf{1}_A$ is a linear space. 
In fact, it is not difficult to see that $\overline{\mbf{T}_n(A)}- \mbf{1}_A$ is isomorphic to the Lie algebra of the maximal torus $\mbf{T}_n(A)$ 
of $\mbf{G}(P,A)$.

Since we can view $\Inc{P,A}-\mbf{1}_A$ as a subset of $\text{End}(\C^n)$, we can apply the additive Jordan decomposition,~\cite[Proposition 2.4.4]{Springer}, to its elements. 
For $y\in \Inc{P,A}$, let $(y_s- \mbf{1}_A)+y_n$ denote the additive Jordan decomposition of $y- \mbf{1}_A$, where $y_s$ is a diagonal matrix and $y_n$ is a strictly upper triangular (nilpotent) matrix. 
We notice that $y_n$ is an element of the Lie algebra of $\mbf{U}_n(P)$.
Now we define the following map:
\begin{align*}
f : \Inc{P,A} &\longrightarrow  \text{Lie}(\mbf{T}_n(A))\oplus  \text{Lie}(\mbf{U}_n(P)) \\
y &\longmapsto (y_s-\mbf{1}_A) + y_n.
\end{align*}
It is easy to check that $f$ is a birational map of varieties. 
Furthermore, it follows from the uniqueness of the additive Jordan decomposition that $f$ is bijective. 
Since the affine space $\text{Lie}(\mbf{T}_n(A))\oplus  \text{Lie}(\mbf{U}_n(P))$ is a normal variety, 
by the Zariski's Main Theorem~\cite[Theorem 5.2.8]{Springer}, we see that $f$ is an isomorphism.
This finishes the proof of our assertion. 
\end{proof}

\section{Bipartite Posets}\label{S:Bipartite}

Let $P$ be a ranked poset such that every interval of $P$ has at most two elements.
If, in addition, $P$ is connected, then we will call $P$ a {\em bipartite poset}.
This terminology is justified by the fact that the Hasse diagram of $P$ is a directed bipartite graph.

The {\em type} of a bipartite poset $P$ is the pair $(k,m)\in \Z_+\times \Z_+$, where $k$ is the number of minimal elements of $P$ and $m$ is the number of maximal elements of $P$. 
A {\em complete bipartite poset of type $(k,m)$} is a bipartite poset of type $(k,m)$ such that every minimal element of $P$ is covered by all maximal elements of $P$.

\begin{Example}\label{E:Mposet}
Let us consider the bipartite posets in Figure~\ref{F:Mposet}. 
The poset $P$ on the left is a bipartite poset of type $(3,2)$. 
The poset $Q$ on the right is a complete bipartite poset of type $(3,2)$. 
\begin{figure}[htp]
\centering
\begin{tikzpicture}[scale=.8]
\begin{scope} [xshift=-5cm]

\node at (-2,0) {$\bullet$};
\node at (0,0){$\bullet$}; 
\node at (2,0){$\bullet$};
\node at (-1,1) {$\bullet$};
\node at (1,1) {$\bullet$};
\draw[-, thick] (-2,0) to  (-1,1);
\draw[-, thick] (2,0) to  (1,1);
\draw[-, thick] (0,0) to  (-1,1);
\draw[-, thick] (0,0) to  (1,1);

\node at (-2,-.5) {$x_1$};
\node at (0,-.5) {$x_2$};
\node at (2,-.5) {$x_3$};
\node at (-1,1.5) {$x_4$};
\node at (1,1.5) {$x_5$};
\node at (0,-1.5) {$P$};

\end{scope}

\begin{scope} [xshift=5cm]

\node at (-2,0) {$\bullet$};
\node at (0,0){$\bullet$}; 
\node at (2,0){$\bullet$};
\node at (-1,1) {$\bullet$};
\node at (1,1) {$\bullet$};
\draw[-, thick] (-2,0) to  (-1,1);
\draw[-, thick] (2,0) to  (1,1);
\draw[-, thick] (0,0) to  (-1,1);
\draw[-, thick] (0,0) to  (1,1);
\draw[-, thick] (-2,0) to  (1,1);
\draw[-, thick] (2,0) to  (-1,1);
\node at (-2,-.5) {$x_1$};
\node at (0,-.5) {$x_2$};
\node at (2,-.5) {$x_3$};
\node at (-1,1.5) {$x_4$};
\node at (1,1.5) {$x_5$};

\node at (0,-1.5) {$Q$};
\end{scope}

\end{tikzpicture}
\caption{Two bipartite posets.}
\label{F:Mposet}
\end{figure}
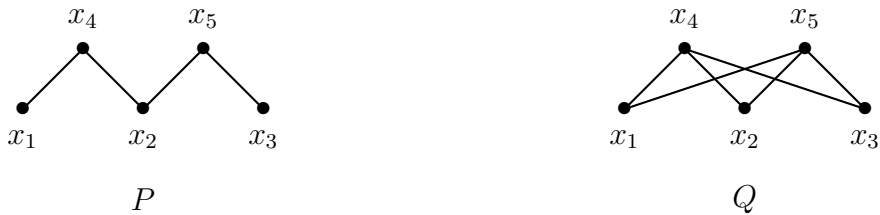

A generic element of $\Inc{P}$ is of the form
\begin{align*}
\begin{bmatrix}
a_{11} & 0 & 0 & a_{14} & 0 \\
0& a_{22} & 0 & a_{24}& a_{25} \\
0& 0 &  a_{33}& 0& a_{35} \\
0& 0 & 0  & a_{44}& 0 \\
0 & 0 & 0 & 0 & a_{55}
\end{bmatrix} \qquad \text{for some $a_{i,j} \in \C$.}
\end{align*}
Similarly, the elements of the incidence monoid of $Q$ are of the form 
\begin{align}\label{A:Mposet}
\begin{bmatrix}
b_{11} & 0 & 0 & b_{14} & b_{15} \\
0& b_{22} & 0 & b_{24}& b_{25} \\
0& 0 &  b_{33}& b_{3,4}& b_{35} \\
0& 0 & 0  & b_{44}& 0 \\
0 & 0 & 0 & 0 & b_{55}
\end{bmatrix}, \qquad \text{where $b_{i,j} \in \C$.}
\end{align}
We notice here that $\Inc{P}$ is a submonoid of $\Inc{Q}$.
\end{Example}

\begin{Remark}\label{R:closedsubmonoid}
The observation we have at the end of this example holds in a greater generality. 
Indeed, it is easy to verify that 
if $P$ is a bipartite poset of type $(k,m)$, and $Q$ is the complete bipartite poset of type $(k,m)$, 
then $\Inc{P}$ is an algebraic submonoid of $\Inc{Q}$. 
\end{Remark}

Recall that a semigroup $S$ is called an {\em orthodox} semigroup if $E(S)$ is a subsemigroup of $S$. 
Finding examples of posets such that $\Inc{P}$ is not an orthodox semigroup is not difficult. 
Nonetheless, we are able to show that all antichain submonoids of certain posets are orthodox semigroups.

\begin{Theorem}\label{T:bipartiteorthodox}
Let $P$ be a ranked poset such that every interval in $P$ has at most two elements. 
Let $A$ be an antichain in $P$. Then the incidence monoid $\Inc{P,A}$ is an orthodox semigroup. 
\end{Theorem}

\begin{proof}
Since $P$ is a disjoint union of bipartite posets, 
in light of Proposition~\ref{P:connectedcomponents}, it suffices to prove our claim for bipartite posets. 
We proceed with the assumption that $P$ is a bipartite poset of type $(k,m)$.
Then the elements of the incidence monoid $\Inc{P}$ are of the form 
\begin{align}\label{A:PorQ}
X:=
\begin{bmatrix}
D_1 & B \\
0 & D_2
\end{bmatrix},
\end{align}
where $D_1$ (resp. $D_2$) is a $k\times k$ diagonal matrix (resp. $m\times m$ diagonal matrix),
and $B$ is a $k\times m$ matrix. 
Clearly, $X$ is invertible if and only if both of $D_1$ and $D_2$ are invertible matrices. 
Since $\mbf{G}(P)$ is a subgroup of $\mathbf{B}_n$, the unipotent radical of $\mbf{G}(P)$ is given by 
the intersection $\mbf{G}(P)\cap R_u(\mbf{B}_n)$. It is easy to check that $R_u(\mbf{B}_n)=\mbf{U}_n$.
By using~\cite[Proposition 2.4.12]{Springer}, we see that $R_u(\mbf{G}(P))\subseteq \mbf{U}_n$.
It follows that we have
\begin{align*}
R_u(\mbf{G}(P)) =
\left\{
\begin{bmatrix}
\mbf{1}_k & B \\
0 & \mbf{1}_m
\end{bmatrix} :\
B\in \mbf{Mat}_{k,m}
\right\}.
\end{align*}
Notice that 
\begin{align*}
\begin{bmatrix}
\mbf{1}_k & B \\
0 & \mbf{1}_m
\end{bmatrix}
\begin{bmatrix}
\mbf{1}_k & B' \\
0 & \mbf{1}_m
\end{bmatrix}
=
\begin{bmatrix}
\mbf{1}_k & B+B' \\
0 & \mbf{1}_m
\end{bmatrix}
\end{align*}
for every $B$ and $B'$ from $\mbf{Mat}_{k,m}$.
It follows that $R_u(\mbf{G}(P))$ is an abelian group. 

By using an argument as in the previous paragraph, 
we see that the unipotent radical of $\mbf{G}(P,A)$ is a subgroup of $\mbf{U}_n \cap \mbf{G}(P)$.
Since $R_u(\mbf{G}(P))$ is the maximal connected normal unipotent subgroup of $\mbf{G}(P)$, 
we see that the unipotent radical of $\mbf{G}(P,A)$ is contained in $R_u(\mbf{G}(P))$.
In particular, we see that $R_u(\mbf{G}(P,A))$ is an abelian group as well. 
It follows from Lemma~\ref{L:smooth} that $\Inc{P,A}$ is a smooth, hence normal, variety. 
The rest of the proof now follows from a result of Renner,~\cite[Theorem 3.2 (c)]{Renner1989}.
\end{proof}

\subsection{The idempotent monoid of a complete bipartite poset.}

Since it is a rather important special case, in this subsection, we determine completely the structure of the irreducible components of the idempotent monoid of a complete bipartite poset $Q$ of type $(k,m)$.
We set, as usual, $n:=k+m$.

\begin{Lemma} \label{L:toporbottom}
Let $X=(X_{ij})_{i,j=1,\dots, n}$ be an idempotent in the incidence monoid $\Inc{Q}$. 
If $X_{ij} \neq 0$ for some $1\leq i < j \leq n$, then $\{ X_{ii},X_{jj} \} = \{ 0,1\}$. 
\end{Lemma}
\begin{proof}
Notice that $X_{ij}\neq 0$ for some $1\leq i < j \leq n$, then $i \leq k < j$. 
Then the $(i,j)$-th entry of $X^2$ is given by $X_{ii}X_{ij} + X_{ij}X_{jj}$.
Since $X^2 = X$, we have $X_{ii}X_{ij} + X_{ij}X_{jj}= X_{ij}$.
Equivalently, we have $X_{ii}+X_{jj} = 1$. 
The proof now follows from the fact that $\{ X_{ii},X_{jj}\}\subseteq \{ 0 ,1\}$.
\end{proof}

The proof of the following corollary follows easily from Lemma~\ref{L:toporbottom}.
\begin{Corollary}\label{C:toporbottom}
Let $X=(X_{ij})_{i,j=1,\dots, n}$ be an idempotent in the incidence monoid $\Inc{Q}$. 
If $X_{ii}=X_{jj}$ for some $i$ and $j$ such that $1\leq i \leq k <j \leq n$, then we have $X_{ij}=0$. 
\end{Corollary}

Next, in the case of a complete bipartite poset $Q$, 
we give a combinatorial formula for the dimension of the irreducible components of $E(\Inc{Q})$.

\begin{Proposition}\label{P:thirdmain}
Let $\ind\subseteq [n]$. 
Then the dimension of the fiber $p^{-1}(\mbf{1}_\ind) \subseteq E(\Inc{Q})$ is given by 
\begin{align*}
\dim p^{-1}(\mbf{1}_\ind) = km - ab - (k-a)(m-b) = kb+am- 2ab,
\end{align*}
where $a= | \ind \cap \{1,\dots, k\}|$ and $b= | \ind \cap \{k+1,\dots, n\}|$.
\end{Proposition}

\begin{proof}
We set 
\begin{align*}
\ind_1:= \{1,\dots, k\} \cap \ind \ \text{ and } \ \ind_2 := \{k+1,\dots, n\} \cap \ind.
\end{align*} 
Let $Y$ be an idempotent in $p^{-1}(\mbf{1}_\ind)$. 
Let $(Y_{ij})_{i,j=1,\dots, n}$ denote the matrix of $Y$. 
Since $p(Y) = \mbf{1}_\ind$, and $Y$ is upper triangular, we have the following restrictions on the `lower triangular' entries of $Y$:
\begin{align*}
Y_{ij} =
\begin{cases}
0 & \ \text{ if $i> j$};\\
0 & \ \text{ if $i=j$ and $i\notin \ind$};\\
1 & \ \text{ if $i=j$ and $i\in \ind$}.
\end{cases}
\end{align*}
The restrictions on the strictly upper triangular entries of $Y$ come from Corollary~\ref{C:toporbottom}:
\begin{align*}
Y_{ij} =
\begin{cases}
0 & \ \text{ if $i\in \ind_1$ and $j\in \ind_2$};\\
0 & \ \text{ if $i< j$, $i\notin \ind_1$, and $j\notin \ind_2$}.
\end{cases}
\end{align*}
There are no other restrictions of the entries of $Y$. 
Therefore, if $i\in \ind_1$, then for each $j\in \{k+1,\dots, n\}\setminus \ind_2$, the entry $Y_{ij}$ can be any element of $\C$. 
Likewise, if $j\in \ind_2$, then for each $i\in \{1,\dots, k\}\setminus \ind_1$, the entry $Y_{ij}$ can be any element of $\C$. 
The zero entries of $Y$ that we listed here will be called the {\em absolutely-zero entries of $Y$}. 
We observe that for each $i\in \ind_1$, the number of zero entries in the $i$-th row of $Y$ are completely determined by the elements $j\in \ind_2$. Similarly, 
for each $i\in \{1,\dots,k\} \setminus \ind_1$, the number of zero entries in the $i$-th row of $Y$ are completely determined by the elements $j\in \{k+1,\dots,n\} \setminus \ind_2$. 
Therefore, the number of (possibly) nonzero entries of $Y$ is given by $km$ minus the total number of absolutely-zero entries of $Y$. 
This finishes the proof of our second assertion. 
\end{proof}

\begin{Example}
Let $\ind$ denote the subset $\ind = \{2,4,6,7\} \subseteq \{1,\dots, 8\}$. 
Let $Q$ denote the complete bipartite poset of type $(3,5)$. 
Then $p^{-1}(\mbf{1}_\ind)$ consists of matrices of the form 
\begin{align*}
\begin{bmatrix}
0 & 0 & 0 & * & 0 & * & * & 0 \\
0 & 1 & 0 & 0 & * & 0 & 0 & * \\
0 & 0 & 0 & * & 0 & * & * & 0 \\
0 & 0 & 0 & 1 & 0 & 0 & 0 & 0 \\
0 & 0 & 0 & 0 & 0 & 0 & 0 & 0 \\
0 & 0 & 0 & 0 & 0 & 1 & 0 & 0 \\
0 & 0 & 0 & 0 & 0 & 0 & 1 & 0 \\
0 & 0 & 0 & 0 & 0 & 0 & 0 & 0 \\
\end{bmatrix},
\end{align*}
where $*$ stands for an arbitrary scalar from $\C$. 
As in the proof of Proposition~\ref{P:thirdmain}, we set $\ind_1:=\ind\cap \{1,2,3\}$ and $\ind_2:=\ind\cap \{4,5,6,7,8\}$.
Since $\ind_1 = \{2\}$ and $\ind_2 = \{4,6,7\}$, we have $a= 1$ and $b=3$. 
Therefore, the dimension of $p^{-1}(\mbf{1}_\ind)$ is given by 
\begin{align*}
km - ab - (k-a)(m-b) = 8.
\end{align*}
\end{Example}
 \medskip

We emphasize once more the fact that for any arbitrary poset $P$ and an antichain $A$ in $P$, 
the idempotent variety of $\Inc{P,A}$ need not be a subsemigroup.
But Proposition~\ref{P:thesameidempotents} together with Theorem~\ref{T:firstmain} show 
that every irreducible component of $E(\Inc{P,A})$ is a subsemigroup. 
In the case of a complete bipartite poset $Q$, and for any antichain $A\subseteq Q$, 
Theorem~\ref{T:secondmain} shows that $\Inc{Q,A}$ is an orthodox  semigroup.
To fully describe the structure of the semigroup $E(\Inc{Q,A})$, we want to answer the following question: 
\medskip
\begin{Question}
How do the subsemigroups $p^{-1}(\mbf{1}_\ind) \subseteq E(\Inc{Q,A})$ fit together?
\end{Question}

First of all, we notice from the homomorphism property of the morphism $p$, for every pair of subsets $\ind,\ind'$ from $[n]$, 
there is an inclusion $p^{-1}(\mbf{1}_\ind) p^{-1}(\mbf{1}_{\ind'})  \subseteq p^{-1}(\mbf{1}_{\ind\cap \ind'})$ 
in $E(\Inc{P})$, where $P$ is any finite poset.

In light of Remark~\ref{R:twoitems}, to answer our question it suffices to work with a maximal antichain $A$ in $Q$.
This means that, if $Q$ is given by $Q=\{x_1,\dots, x_n\}$, then $A$ is either $\{x_1,\dots, x_k\}$ or $\{x_{k+1},\dots, x_n\}$. 
Without losing generality, we will proceed with the assumption that $A=\{x_{k+1},\dots, x_n\}$. 
In this case, a diagonal idempotent $\mbf{1}_\ind$ is a member of $\Inc{Q,A}$ if and only if the inclusions $[k]\subseteq \ind\subseteq [n]$ hold.

\begin{Theorem}\label{T:componentmultiplication}
Let $Q$ be a complete bipartite poset of type $(k,n-k)$ and let $A$ be the maximal antichain $\{x_{k+1},\dots, x_n\}$ in $Q$.
The irreducible components of $E(\Inc{Q,A})$ are given by $p^{-1}(\mbf{1}_\ind)$, where $[k]\subseteq \ind\subseteq [n]$ hold.
Furthermore, for every two such subsets $\ind$ and $\ind'$, we have 
$p^{-1}(\mbf{1}_\ind) p^{-1}(\mbf{1}_{\ind'}) = p^{-1}(\mbf{1}_{\ind\cap \ind'})$.
\end{Theorem}
\begin{proof}
We already proved our first assertion, so, we proceed with the proof of the second assertion. 
Under our assumptions, our claim follows from a direct matrix multiplication.
Let $Y \in p^{-1}(\mbf{1}_\ind)$ and $Z\in p^{-1}(\mbf{1}_{\ind'})$.
Then the columns of $YZ$ are determined as follows.
Let $Y_i$ (resp. $Z_i$) denote the $i$-th column of $Y$ (resp. of $Z$) for $i\in [n]$.
Likewise, we denote by $(YZ)_i$ the $i$-th column of $YZ$ for $i\in [n]$.
Then by a direct computation and using Corollary~\ref{C:toporbottom} we see that, for every $i\in [n]$, the $i$-th column of $YZ$ is given by 
\begin{align}\label{A:determiningcolumns}
(YZ)_i =
\begin{cases}
Z_i & \text{if the $i$-th entry of the column $Z_i$ is zero},\\
Y_i & \text{otherwise}.
\end{cases}
\end{align}
We notice that both of the components $p^{-1}(\mbf{1}_\ind)$ and $p^{-1}(\mbf{1}_{\ind'})$ are affine spaces. 
Indeed, this follows from the fact the matrix entries of the elements of $p^{-1}(\mbf{1}_\ind)$ and $p^{-1}(\mbf{1}_{\ind'})$ can be freely chosen from the underlying field.
The calculation in (\ref{A:determiningcolumns}) shows that the product $p^{-1}(\mbf{1}_\ind)p^{-1}(\mbf{1}_{\ind'})$ is also an affine space of dimension $k (n- |\ind\cap \ind'|)$.
Since $p^{-1}(\mbf{1}_\ind)p^{-1}(\mbf{1}_{\ind'}) \subseteq p^{-1}(\mbf{1}_{\ind\cap \ind'})$, to finish the proof, 
it remains to show that $\dim p^{-1}(\mbf{1}_{\ind\cap \ind'}) = k (n- |\ind\cap \ind'|)$.
We will use Proposition~\ref{P:thirdmain}. 

By our assumptions on $\ind$ and $\ind'$ we have $[k] \subseteq \ind\cap \ind'$. 
Therefore, $a$ and $b$ of Proposition~\ref{P:thirdmain} are given by 
\begin{align*}
a= k\qquad\text{and}\qquad b= |\ind\cap \ind'|-k.
\end{align*}
Then, by the same proposition we have
\begin{align*}
\dim p^{-1}(\mbf{1}_{\ind\cap \ind'}) &= kb +a(n-k) -2ab \\
&= k(|\ind\cap \ind'|-k) + kn -k^2 - 2k(|\ind\cap \ind'|-k) \\
&= kn - k|\ind\cap \ind'|.
\end{align*}
This finishes the proof of our theorem.
\end{proof}

\begin{Example}
Let $Y\in p^{-1}(\mbf{1}_{\{1,2,3,6,7,8\}})$ and $Z\in p^{-1}(\mbf{1}_{\{1,2,3,5,8\}}) $ denote the following idempotent matrices: 
\begin{align*}
Y:= \begin{bmatrix}
1 & 0 & 0 & a_{14} & a_{15} & 0 & 0 & 0 \\
0 & 1 & 0 & a_{24} & a_{25} & 0 & 0 & 0 \\
0 & 0 & 1 & a_{34} & a_{35} & 0 & 0 & 0 \\
0 & 0 & 0 & 0 & 0 & 0 & 0 & 0 \\
0 & 0 & 0 & 0 & 0 & 0 & 0 & 0 \\
0 & 0 & 0 & 0 & 0 & 1 & 0 & 0 \\
0 & 0 & 0 & 0 & 0 & 0 & 1 & 0 \\
0 & 0 & 0 & 0 & 0 & 0 & 0 & 1 
\end{bmatrix}\ \text{ and } \ 
Z:=  \begin{bmatrix}
1 & 0 & 0 & b_{14} & 0 & b_{16} & b_{17} & 0 \\
0 & 1 & 0 & b_{24} & 0 & b_{26} & b_{27} & 0 \\
0 & 0 & 1 & b_{34} & 0 & b_{36} & b_{37} & 0 \\
0 & 0 & 0 & 0 & 0 & 0 & 0 & 0 \\
0 & 0 & 0 & 0 & 1 & 0 & 0 & 0 \\
0 & 0 & 0 & 0 & 0 & 0 & 0 & 0 \\
0 & 0 & 0 & 0 & 0 & 0 & 0 & 0 \\
0 & 0 & 0 & 0 & 0 & 0 & 0 & 1 
\end{bmatrix}.
\end{align*}
Then we have 
$
YZ=
 \begin{bmatrix}
1 & 0 & 0 & b_{14} & a_{15} & b_{16} & b_{17} & 0 \\
0 & 1 & 0 & b_{24} & a_{25} & b_{26} & b_{27} & 0 \\
0 & 0 & 1 & b_{34} & a_{35} & b_{36} & b_{27} & 0 \\
0 & 0 & 0 & 0 & 0 & 0 & 0 & 0 \\
0 & 0 & 0 & 0 & 0 & 0 & 0 & 0 \\
0 & 0 & 0 & 0 & 0 & 0 & 0 & 0 \\
0 & 0 & 0 & 0 & 0 & 0 & 0 & 0 \\
0 & 0 & 0 & 0 & 0 & 0 & 0 & 1 
\end{bmatrix}\in p^{-1}(\mbf{1}_{\{ 1,2,3,8\}}).
$
\end{Example}

\section{The Cross-Section Lattice}\label{S:Cross}

Recall from Subsection~\ref{SS:Preliminaries:Semigroups} that a regular irreducible linear algebraic monoid $M$ has a cross-section lattice. 
Our goal in this section is to determine a cross-section lattice of $\Inc{P,A}$, where $P$ is a bipartite poset, and $A$ is an antichain.
As our earlier Example~\ref{E:Mposet} indicates, by adding cover relations to our bipartite poset, we enlarge the underlying incidence monoid without changing its maximal torus, hence, its diagonal idempotents. At the same time, our earlier Remark~\ref{R:twoitems}
shows that an enlargement of an antichain adds more idempotents to the underlying antichain monoid.
Since there is no information lost when we work with complete bipartite posets and their maximal antichains, we proceed with this assumption.

We fix the following notation: $Q=\{x_1,\dots, x_n\}$ is a complete bipartite poset of type $(k,m)$, so, $n=m+k$. 
The letter $A$ denotes the maximal antichain that consists of maximal elements, $A:=\{x_{k+1},\dots, x_n\}$. 
If there is no danger for confusion, we will denote the antichain monoid $\Inc{Q,A}$ by $M$; the unit group of $M$ will be denoted by $G$.

\begin{Theorem}\label{T:Cross-section}
In the above notation, the cross-section lattice $\Lambda$ of $M$ is isomorphic to the Boolean lattice of all subsets of $[m]$.
\end{Theorem}
\begin{proof} 
Since $M$ is an irreducible linear algebraic monoid, for $x\in M$, the $\gj$-class of $x$ is given by $GxG$.
In particular, for every idempotent $e\in E(M)$, we have $J_e = GeG$. 
At the same time, two idempotents of $M$ are in the same $\gj$-class if and only if they are $G$-conjugate. 
By Proposition~\ref{P:thesameidempotents}, we know that the diagonal idempotents $\mbf{1}_\ind$, where $[k]\subseteq \ind\subseteq [n]$,
are the representatives of the $G$-conjugacy classes in $E(M)$.

Since $M$ is an irreducible regular monoid, it has a cross-section lattice $\Lambda$. 
Furthermore, all cross-section lattices are conjugate to each other. 
Without loss of generality we will work with the cross-section lattice that is contained in the maximal torus $\overline{\mbf{T}_n(A)}$.

Recall that $E(\overline{\mbf{T}_n(A)})=\{\mbf{1}_\ind : [k]\subseteq \ind\subseteq [n]\} $.
Now, since every element of $\Lambda$ intersects a unique $\gj$-class and every $\gj$-class is represented by a unique diagonal idempotent $\mbf{1}_\ind$ such that $[k]\subseteq \ind\subseteq [n]$, we see that the elements of $\Lambda$ are in bijection with the subsets of $[n-k]=[m]$. 
Finally, it is easy to check that the natural order on the idempotents in 
$\{\mbf{1}_\ind:\ [k]\subseteq \ind\subseteq [n]\}$ is the one that corresponds to the inclusion order on their indices.
Hence we recover $\Lambda$ as the Boolean lattice $([m],\subseteq)$. 
This finishes the proof of our theorem. 
\end{proof}

The proof of our theorem shows that a cross-section lattice for an antichain monoid can be chosen inside the closure of the diagonal torus. We isolate this fact as a corollary to refer to it later. 

\begin{Corollary}\label{C:Cross-section}
In the above notation, the set
\[\Lambda=\left\{ \mbf{1}_\ind : [k]\subseteq \ind\subseteq [n] \right\} \]
is a  cross-section lattice of $M$.
\end{Corollary}
\begin{proof}
The proof is contained in the second paragraph of the proof of Theorem~\ref{T:Cross-section}.
\end{proof}

\begin{Problem}
Is it true that the weak cross-section lattice of an incidence monoid is a cross-section lattice if and only if it is a regular monoid? 
\end{Problem}

\subsection{The $\gh$-classes.}

In the remainder of this section, we will determine the structures of the $\gh$-classes of (the elements of the cross-section lattice of) 
$\Inc{Q,A}$, where $Q=\{x_1,\dots, x_n\}$ is a complete bipartite poset of type $(k,m)$, and $A=\{x_{k+1},\dots, x_n\}\subseteq Q$ is the maximal antichain. 
We begin with some general remarks.

Let $M$ be an incidence monoid. 
Let $e$ be an idempotent from $E(M)$. 
Since $G$ is solvable, the Weyl group $W$ of $G$ is trivial. 
By~\cite[Corollary 6.34]{Putcha}, we know that the following are equivalent:
\begin{enumerate}
\item $1= w(e) := | \{ w e w^{-1} : w\in W \}|$. 
\item $G=C_G^l (e) C_G^r(e)$, where $C_G^l(e)$ (resp. $C_G^r(e)$) is the left (resp. right) centralizer of $e$ in $G$.
\item $eGe$ is the $\gh$-class of $e$.
\end{enumerate}
Let $f$ be another idempotent that is $G$-conjugate to $e$.
Let $g\in G$ be such that $f=geg^{-1}$.
Then we have $fGf = g e g^{-1} G g eg^{-1} = g eGe g^{-1}$.
It is easily seen that the map $x\mapsto gxg^{-1}$, where $x\in eGe$, defines an isomorphism between $eGe$ and $fGf$.
From these observations we conclude that to calculate the $\gh$-class of an idempotent of $M$, it suffices to focus on the groups $eGe$, where $e$ is from the cross-section lattice of $M$. 
\medskip

We now specialize to the antichain monoid $M:=\Inc{Q,A}$. 
Since $M$ is completely regular, by~\cite[Proposition 4.1.1]{Howie}, every $\gh$-class in $M$ is a group.
Hence, every $\gh$-class contains an idempotent of $M$. 
By our observations in the previous paragraph, we see that any $\gh$-class in $M$ is isomorphic to one of the 
$\gh$-classes of the form $eGe$, where $e \in \Lambda$. 
By Corollary~\ref{C:Cross-section}, we know that 
\begin{align*}
\Lambda = \{ \mbf{1}_\ind :\ [k]\subseteq \ind \subseteq [n] \}.
\end{align*}
\begin{Notation}
Let $\mbf{1}_\ind$ be an idempotent from $\Lambda$.
If $\mbf{1}_\ind$ is given by $\mbf{1}_\ind = \text{diag}(a_1,\dots, a_n)$, then we set 
\begin{align*}
\widetilde{\mbf{1}_\ind}:=\text{diag}(a_{k+1},\dots, a_n).
\end{align*}
We denote by $\TM_\ind$ the set of diagonal matrices $\diag \in \TM_m$  such that $\diag =\diag'\widetilde{\mbf{1}_\ind}$ for some $\diag '\in \mbf{T}_m$.
\end{Notation}
The elements of the unit group $G$ are the invertible matrices of the form 
$Y: = \begin{bmatrix}\mbf{1}_k &  B  \\ 0 & \diag \end{bmatrix}$,
where $\diag $ is an invertible diagonal $m\times m$ matrix, and $B\in \mbf{Mat}_{k,m}$.

\begin{Proposition}\label{P:Hclass}
Let $\mbf{1}_\ind$ be an idempotent from $\Lambda$ for some $\ind$ such that $[k]\subseteq \ind\subset [n]$.
Then the $\gh$-class of $\mbf{1}_\ind$ is given by 
\begin{align*}
\mbf{1}_\ind G \mbf{1}_\ind= 
\left\{
\begin{bmatrix}
\mbf{1}_k & B \\
0 & \diag 
\end{bmatrix}:  
\begin{matrix} 
\text{$B\in \mbf{Mat}_{k,m}$ is such that $B\widetilde{\mbf{1}_\ind}=B$} \\
\text{$\diag \in \TM_\ind$}
\end{matrix}
\right\}.
\end{align*}
\end{Proposition}
\begin{proof}
This follows from a direct computation by matrix multiplication. We omit the details. 
\end{proof}


Proposition~\ref{P:Hclass} suggests that some of the $\gh$-classes in different $\gd$-classes can be isomorphic to each other as well. 
We make this observation more precise by our next result.

\begin{Theorem}\label{T:isomorphicHclasses}
Let $\mbf{1}_\ind$ be an idempotent from $\Lambda$ for some $\ind$ such that $[k]\subseteq \ind\subseteq [n]$.
Let $a:= |\ind|$. 
If $a>k$, then the $\gh$-class of $\mbf{1}_\ind$ is isomorphic to the unit group $\mbf{G}(Q(a),A(a))$,
where $Q(a)$ is the complete bipartite poset $Q(a):=\{x_1,\dots, x_a\}$, and $A(a)$ is the maximal antichain 
$\{x_{k+1},\dots, x_a\}$. 
\end{Theorem}

\begin{proof}
Let $G(a)$ denote the unit group of $\mbf{G}(Q(a),A(a))$.
Strictly speaking, $G(a)$ is contained in the standard Borel subgroup $\mathbf{B}_a$. 
We consider the map defined by 
\begin{align*}
\psi: \mbf{1}_\ind G \mbf{1}_\ind  &\longrightarrow  G(a) \\
\begin{bmatrix}
\mbf{1}_k & B  \\
0 & \diag   
\end{bmatrix} &\longmapsto 
\begin{bmatrix}
\mbf{1}_k & B' \\
0 & \diag' 
\end{bmatrix},
\end{align*}
where $B'\in \mbf{Mat}_{k, n-a}$ is the matrix obtained from $B$ 
by deleting its columns indexed by the elements of $[n]\setminus \ind$. 
Likewise, $\diag'\in \mathbf{T}_{n-a}$ is the matrix obtained from $\diag$ by deleting its columns as well as rows that correspond to the elements of $[n]\setminus \ind$. 
Since $\psi$ is defined by deleting only the zero columns (and zero rows) it is an isomorphism of varieties. 
We claim that $\psi$ is a group homomorphism. 
To see this, let 
\begin{align*}
X:=
\begin{bmatrix}
\mbf{1}_k & B_1 \\
0 & \diag_1 
\end{bmatrix}\ \text{ and } \
Y:=
\begin{bmatrix}
\mbf{1}_k & B_2 \\
0 & \diag_2
\end{bmatrix}
\end{align*}
be two elements from $\mbf{1}_\ind G \mbf{1}_\ind$. 
Then 
$
XY:=
\begin{bmatrix}
\mbf{1}_k & B_2+B_1T_2 \\
0 & \diag_1\diag_2 
\end{bmatrix}.
$
It is easy to check that $(B_2+B_1\diag_2)' = B_2' +B_1'\diag_2'$ and $(\diag_1\diag_2)' = \diag_1'\diag_2'$. 
In other words, we have $\psi (XY) = \psi (X)\psi(Y)$. 
This finishes the proof of our assertion.
\end{proof}

\begin{Corollary}\label{C:fourthmain}
Let $\mbf{1}_\ind$ and $\mbf{1}_{\ind'}$ be two elements from $\Lambda$.
If $|\ind| = |\ind'| > k$, then the $\gh$-classes of $\mbf{1}_\ind$ and $\mbf{1}_{\ind'}$ are isomorphic. 
Furthermore, there are exactly $m+1$ nonisomorphic $\gh$-classes in $M$. 
\end{Corollary}
\begin{proof}
Let $G(a)$ denote the unit group $\mbf{G}(Q(a),A(a))$, where $a=|\ind|$, and $Q(a)$ is as defined in Theorem~\ref{T:isomorphicHclasses}.
Then we know that both of the $\gh$-classes $\mbf{1}_\ind G \mbf{1}_\ind$ and 
$\mbf{1}_{\ind'} G \mbf{1}_{\ind'}$ are isomorphic to $G(a)$. 
The spectrum of cardinalities $|\ind|$, where $\mbf{1}_\ind \in \Lambda$, is given by the set $\{k,k+1,\dots, n\}$. 
By comparing their maximal tori, we see immediately that $G(a) \ncong G(b)$ if $k\leq a\neq b \leq n$.  
\end{proof}

\section{Green's Relations Revisted}\label{S:Green}

In the previous sections we determined the $\gh$-classes (of the elements of the cross-section lattices) 
as well as the idempotent semigroups of the antichain monoids of (complete) bipartite posets.
In this section, we will characterize the remaining Green's relations. 
As before, let $Q=\{x_1,\dots, x_n\}$ be a complete bipartite poset of type $(k,m)$.
Let $A=\{x_{k+1},\dots, x_n\}$ be a maximal antichain in $Q$.
We proceed to determine the inverses of the elements of $M:=\Inc{Q,A}$. Recall that $M$ is a completely regular semigroup \cite[Proposition~2.4]{Can2021}.

Let $X$ be an element of $M$.
Then $X$ is of the form 
$X=\begin{bmatrix}
\mathbf{1}_k &  B \\ 
0 & \diag
\end{bmatrix}$,
for some $B\in \mbf{Mat}_{k,m}$ and  $\diag \in\TM_L$, with $\ind$ such that $[k]\subseteq \ind\subseteq [n]$. 
Since $\TM_m$ is a commutative semigroup, there is a unique inverse $\diag'\in V(\diag)$ such that 
$\diag\diag'=\diag'\diag=\widetilde{\mbf{1}_\ind}$.
After using the defining equations of $V(X)$, we find that 
\begin{align*}
V(X)=
\left \{
\begin{bmatrix}
\mathbf{1}_k &  B' -B\diag'\\ 0 & \diag'
\end{bmatrix}\,:\, B'\widetilde{\mathbf{1}_{\ind}}=0,\,\diag\diag'=\diag'\diag=\widetilde{\mathbf{1}_{\ind}}\right\}.
\end{align*}
Hereafter, if there is no danger for confusion, we denote by $X'$ any element of $V(X)$, and we denote by $X^{-1}$ the unique element in $V(X)$ satisfying (\ref{A:inverseincrs}). With some calculations we obtain 
\begin{align*}
X^{-1}=
\begin{bmatrix}
\mathbf{1}_k &   B(\mathbf{1}_m-\widetilde{\mathbf{1}_{\ind}})-B\diag' \\ 0 & \diag'
\end{bmatrix},\ \textrm{ where }\ \diag\diag'=\diag'\diag=\widetilde{\mathbf{1}_{\ind}}
\end{align*}
and 
\begin{align}\label{A:idempotent}
XX^{-1}=X^{-1}X=\begin{bmatrix} \mbf{1}_k & B(\mathbf{1}_m-\widetilde{\mathbf{1}_{\ind}}) \\ 0 & \widetilde{\mathbf{1}_{\ind}} \end{bmatrix}.
\end{align}

Now it is easy to verify the following characterization of the idempotents of $M$:
\begin{align*}
E(M)=\left\{
\begin{bmatrix}
\mathbf{1}_k &  B \\ 0 & \widetilde{\mathbf{1}_{\ind}}
\end{bmatrix}: \, B\widetilde{\mathbf{1}_{\ind}}=0,\,[k]\subseteq \ind\subseteq [n]\right\}.
\end{align*}

\begin{Remark}\label{R:RightZeroBand}
Let $E\in p^{-1}(\mbf{1}_\ind) \subseteq E(M)$, for some $[k]\subseteq \ind\subseteq [n]$.
Then we have 
\begin{align*}
V(E)=
\left\{
\begin{bmatrix}
\mbf{1}_k & B \\
0 & \widetilde{\mathbf{1}_{\ind}}
\end{bmatrix}\ :\
B \widetilde{\mathbf{1}_{\ind}}=0
\right\}=p^{-1}(\mbf{1}_\ind)=\mbf{G}(Q,A)\cdot \mbf{1}_\ind =J_{\mbf{1}_\ind}\cap E(M),
\end{align*}
by Proposition~\ref{P:thesameidempotents}.
It is easy to check that $V(E)$ is a subsemigroup of $E(M)$, in accordance with the mentioned result. 
Furthermore, $V(E)$ is a right-zero band.
Indeed if $F,G\in V(E)$ then $FG=G$.
\end{Remark}

Any element $X\in M$ is in the same $\gh$-class of the respective idempotent $XX^{-1}=X^{-1}X$. So to verify if any two elements in $M$ are $\gh$-related (resp., $\gr$-, $\gl$- and $\gj$-) is equivalent to verify if the respective idempotents are in the same $\gh$-class (resp., $\gr$-, $\gl$- and $\gj$-). Recall that $\gd=\gj$. So, by the previous remark, $X\gj Y$ if and only if there exists some $[k]\subseteq \ind\subseteq [n]$ such that $XX^{-1}, YY^{-1}\in p^{-1}(\mbf{1}_\ind)$.

Since for any $L$ with $[k]\subseteq \ind\subseteq [n]$, $p^{-1}(\mbf{1}_\ind)$ is a right-zero band, any two idempotents in $p^{-1}(\mbf{1}_\ind)$ are $\gl$-related. Therefore, every $\gj$-class has at most one $\gl$-class. Hence,  $\gl=\gh$ and $\gr=\gj=\gd$.

We summarize our findings in the form of a theorem.

\begin{Theorem}\label{T:summaryform}
Let $X=\begin{bmatrix} \mathbf{1}_k &  B_X \\ 0 & T_X \end{bmatrix}$ be an element of $M$, with $T_X\in \TM_\ind$ for some $L$  such that $[k]\subseteq \ind\subseteq [n]$.
Then we have 
\begin{enumerate}
\item[(1)] $J_X=D_X=R_X=\left\{\begin{bmatrix} \mathbf{1}_k &  B \\ 0 & \diag \end{bmatrix} \in M:\,   \diag \in\TM_\ind  \right\}=J_{\mathbf{1}_{\ind}}$;
\item[(2)] $L_X=H_X=\left\{\begin{bmatrix} \mathbf{1}_k &  B \\ 0 & \diag \end{bmatrix}\in M\ :\ B(\mathbf{1}_{m}-\widetilde{\mathbf{1}_{\ind}})=B_X(\mathbf{1}_{m}-\widetilde{\mathbf{1}_{\ind}}) \text{ and }  \diag \in\TM_\ind \right\}$.
\item[(3)] $J_X$ is a right group. 
\end{enumerate}
\end{Theorem}

\begin{proof}
\noindent (1) From the previous observations we have that $J_X=D_X=R_X$. 
Because $XX^{-1}$ is an idempotent it belongs to $p^{-1}(\mathbf{1}_{\ind})$, for some $L$  with  $[k]\subseteq L\subseteq [n ]$, and so $XX^{-1}\gj \mathbf{1}_{\ind}$. Since $X\gh XX^{-1}$ we get $J_X=J_{\mathbf{1}_{\ind}}$. 

Let $K$ denote the set  $\left\{\begin{bmatrix} \mathbf{1}_k &  B \\ 0 & \diag \end{bmatrix} \in M:\,  \diag \in\TM_\ind  \right\}$. Given $Y\in K$, from equality (\ref{A:idempotent}), we deduce that $YY^{-1}\in p^{-1}(\mathbf{1}_{\ind})\subseteq J_{\mathbf{1}_{\ind}}$. Therefore, $Y\in J_{\mathbf{1}_{\ind}}$ since $Y\gh YY^{-1}$.
Now, if $Y=\begin{bmatrix} \mathbf{1}_k &  B \\ 0 & \diag \end{bmatrix}\in J_{\mathbf{1}_{\ind}}$, then  $YY^{-1}\in p^{-1}(\mathbf{1}_{\ind})$, and so  $\diag \in\TM_\ind$. Therefore, $Y\in K$.   

\noindent (2) From the previous observations we know that $L_X=H_X$.
Let $Y=\begin{bmatrix} \mathbf{1}_k &  B \\ 0 & \diag \end{bmatrix}\in M$. We have $X\gh Y$ if and only if $XX^{-1}=YY^{-1}$, because $H_X$ is a group.
Attending to (\ref{A:idempotent}), we deduce that $Y$ has the form given in the statement.

\noindent (3) We already mentioned the result of Putcha (\cite[Corollary 3.20]{Putcha}) that every $\gj$-class in a closed subsemigroup of $\BM_n$ 
is a completely simple semigroup. By \cite[Proposition~2.4.3]{Lallement}, since $J_X$ is a regular subsemigroup of $M$, then $J_X$ has also only one $\gl$-class. It is known that a semigroup is a right group if and only if it is a completely simple and contains one $\gl$-class \cite[Theorem~1.4.9]{Cain}. Therefore $J_X$ is a right group.
\end{proof}

\begin{Remark}
Let $X$ be an element of the $\gh$-class of an idempotent $E\in E(M)$. 
It is not difficult to verify that the description of $H_X$ in Theorem~\ref{T:summaryform} (2) is in agreement with the elements of the $\gh$-class of $E$, that is, $H_E=EGE$.
\end{Remark}

\begin{Proposition}\label{P:isomorphismJclass}
Let $X$ be an element of $M$ and let $L$, with $[k]\subseteq \ind\subseteq [n]$, be such that $X\in J_{\mathbf{1}_{\ind}}$. The mapping $\rho: H_X\longrightarrow H_{\mathbf{1}_{\ind} }$, given by $Y\mapsto Y\mathbf{1}_{\ind} $, is an isomorphism, and its inverse map is given by $Y\mapsto YXX^{-1}$.
\end{Proposition}
\begin{proof}
The result follows from the fact that $X\gr \mathbf{1}_{\ind} $ and the proof in \cite[Proposition~2.3.6]{Howie}.
\end{proof}

\begin{Remark}
By \cite[pg. 61]{Howie}, given $E \in E(M)$, $J_EE$ is a group and $\phi:J_EE\times E(J_E)\longrightarrow J_E$, $(X,F)\mapsto XF$ is an isomorphism.
\end{Remark}

\section{Conjugacy Relations}\label{S:Conjugacy}

There is a continuous interest on the study of various conjugacy actions for reasons notably rooted in representation theory. 
For reductive monoids, the ordinary conjugation action of the unit group is described by Putcha in~\cite{Put2005}.
For regular semigroups as well as transformation semigroups, and for some specific conjugacy relations that we will discuss below,
the groundbreaking work is done by Kudryavtseva and Mazorchuk in the papers~\cite{Kud06, KM2007, KM2008}.
Our main goal in this section is to compare the conjugation action by the unit group on the monoid and the so called primarily relation $\simp$ credited to Lyndon and Sch\"utzenberger \cite{LS62} by Lallement.
This relation has been considered as one of many other possible generalizations of conjugacy to semigroups - see \cite{AKKM} for different possible notions of conjugacy in semigroups. 

For elements $a$ and $b$ from a semigroup $S$, we say that $a$ and $b$ are {\em primarily related}, or that, $a$ and $b$ are {\em $p$-conjugate}, and write $a\simp b$, if there are $z,w \in S^1$ such that 
\begin{align}\label{A:pconjugacy}
 a=zw ,\, b=wz.
\end{align}
For a completely regular semigroup $S$, it is known that the $p$-conjugacy is transitive and hence an equivalence relation \cite[Corollary 4]{Kud06}.
For this class, we also know that both generalizations of conjugacy to semigroups $\simn$ and $\simtr$ are equal to $\simp$ by \cite[Theorem~6.5]{ABKKMM}  and \cite[Corollary~4.6]{AKKM}.

For our purposes here, it is convenient to work with the following definition. 
Let $a,b\in S$. Then we define 
\begin{align}\label{A:nconjugacy}
a\simn b \Leftrightarrow \exists g,\, h\in S^1:\, ag = gb,\ bh = ha,\ hag = b,\ \textrm{and}\ gbh = a.
\end{align}
Note that if $a\simn b$, then $a\gj b$. 
For the sake of completeness, we include a characterization of the trace conjugacy in the class of completely regular semigroups (see~\cite[Theorem~4.5]{AKKM}):
\begin{eqnarray*}
	a\simtr  b	&\Leftrightarrow& \exists z,\, w\in S^1:\, wzw=w,\,waz=b,\,zbw=a.
\end{eqnarray*}

In the rest of this section, we follow the notation that is setup in the previous section.
In particular, $M$ will denote $\Inc{Q,A}$, where $Q$ and $A$ are as defined at the beginning of Section~\ref{S:Green}. 

By~\cite[Proposition 2.4]{Can2021}, $M$ is a completely regular monoid and so the conjugacy notions $\simn$, $\simtr$ and $\simp$ coincide in $M$. Hence, $p$-conjugate elements must be in the same $\gj$-class. 
We have the following result, which is stated as Theorem~\ref{T:lastannounced} in the introduction. 
\begin{Theorem}\label{T:modifiedlastannounced}
Let $X$ and $Y$ be two elements of $M$.
Then $X\simp Y$ if and only if $\{X,Y\}\subseteq J_{\mathbf{1}_{\ind}}$ for some $L$ with $[k]\subseteq L\subseteq [n]$, and  $X\mathbf{1}_{\ind}$ and  $Y\mathbf{1}_{\ind}$ are conjugate in the group $H_{\mathbf{1}_{\ind}}$.
\end{Theorem}

\begin{proof}
Assume that  $X\simp Y$. Then we already observed that $X\gj Y$. 
So, there exists a set $L$ such that $[k]\subseteq L\subseteq [n]$ and $X,Y\in J_{\mathbf{1}_{\ind}}$. 
By Proposition~\ref{P:isomorphismJclass}, we have $X\mathbf{1}_{\ind},Y\mathbf{1}_{\ind}\in H_{\mathbf{1}_{\ind}}$.
Now, let $Z,W\in M$ be such that $X=ZW$ and $Y=WZ$. 
We notice that 
\[
J_{\mathbf{1}_{\ind}}=J_{X\mathbf{1}_{\ind}}=J_{Y\mathbf{1}_{\ind}}\ \leq \ J_{Z\mathbf{1}_{\ind}},J_{W\mathbf{1}_{\ind}}\ \leq \ J_{\mathbf{1}_{\ind}}\qquad \text{ (see~(\ref{A:partialorderonJ})).}
\]
This means that the elements $Z \mathbf{1}_{\ind}$ and $W\mathbf{1}_{\ind}$ belong to the same $\gj$-class, $J_{\mathbf{1}_{\ind}}$.
Again by Proposition~\ref{P:isomorphismJclass}, we conclude that $Z\mathbf{1}_{\ind},W\mathbf{1}_{\ind}\in H_{\mathbf{1}_{\ind}}$. 
Since $\mathbf{1}_{\ind}$ is the identity element of the group $H_{\mathbf{1}_\ind}$, 
we can write $X\mathbf{1}_{\ind}= \mathbf{1}_{\ind}X\mathbf{1}_{\ind}$,
$Y\mathbf{1}_{\ind}= \mathbf{1}_{\ind}Y\mathbf{1}_{\ind}$, $W\mathbf{1}_{\ind}= \mathbf{1}_{\ind}W\mathbf{1}_{\ind}$,
 and $Z\mathbf{1}_{\ind}=\mathbf{1}_{\ind}Z\mathbf{1}_{\ind}$. 
Therefore, we have 
\[
X\mathbf{1}_{\ind}=(Z W)\mathbf{1}_{\ind}=Z (W\mathbf{1}_{\ind})= Z(\mathbf{1}_\ind W\mathbf{1}_{\ind}) 
= (Z\mathbf{1}_{\ind}) (W\mathbf{1}_{\ind}).
\]
Similarly, we find that $Y\mathbf{1}_{\ind}=W\mathbf{1}_{\ind} Z\mathbf{1}_{\ind}$. 
In other words, the elements $X\mathbf{1}_{\ind}$ and  $Y\mathbf{1}_{\ind}$ are $p$-conjugate in the group $H_{\mathbf{1}_{\ind}}$.
Therefore, they are $H_{\mathbf{1}_{\ind}}$-conjugate.

Conversely, suppose that $\{X,Y\}\subseteq J_{\mathbf{1}_{\ind}}$, for some $L$ with $[k]\subseteq L\subseteq [n]$, and  $X\mathbf{1}_{\ind}$ and  $Y\mathbf{1}_{\ind}$ are conjugate in the group $H_{\mathbf{1}_{\ind}}$.
Let $Z, W\in H_{\mathbf{1}_{\ind}}$ be such that $X\mathbf{1}_{\ind}=ZW$ and $Y\mathbf{1}_{\ind}=WZ$. Consider $Z'=ZYY^{-1}$ and $W'=WXX^{-1}$.
 By  Proposition~\ref{P:isomorphismJclass}, $Z'\in H_Y$ and $W'\in H_X$. Note that $XX^{-1}$ and $YY^{-1}$ are idempotents in the same $\gj$-class being the identities of the groups $H_X$ and $H_Y$, respectively. Attending to Remark~\ref{R:RightZeroBand}, we get
$ Z'W'=ZYY^{-1}WXX^{-1}=ZYY^{-1}XX^{-1}WXX^{-1}=ZXX^{-1}WXX^{-1}=
ZWXX^{-1}=X\mathbf{1}_{\ind}XX^{-1}=XXX^{-1}=X$. Similarly, we have $ W'Z'=WXX^{-1}ZYY^{-1}=WXX^{-1}YY^{-1}ZYY^{-1}=WYY^{-1}ZYY^{-1}=
WZYY^{-1}=Y\mathbf{1}_{\ind}YY^{-1}=YYY^{-1}=Y$.
This finishes the proof of our proposition.
\end{proof}

The previous result shows us that any two elements $X$ and $Y$ of $M$ are $p$-conjugate if and only if they are $\gj$-related (for some $\gj$-class $J_{\mathbf{1}_{\ind}}$) and their images, under the map in Proposition~\ref{P:isomorphismJclass}, are conjugate in the group of units of the idempotent $\mathbf{1}_{\ind}$. Notice that, by Theorem~\ref{T:isomorphicHclasses}, $H_{\mathbf{1}_{\ind}}$ is isomorphic
to the unit group $\mbf{G}(a):=\mbf{G}(Q(a),A(a))$,
where $Q(a)$ is the complete bipartite poset $Q(a):=\{x_1,\dots, x_a\}$, and $A(a)$ is the maximal antichain 
$\{x_{k+1},\dots, x_a\}$. 
As a consequence of our next result, we will show how to identify conjugate elements in such unit groups $\mbf{G}(a)$.
To this end, first, we introduce some terminology. 
Let $G$ be a connected solvable group. 
Let us denote by $T$ a maximal torus in $G$, and by $U$ the unipotent radical of $G$. 
By the multiplicative Jordan decomposition theorem (see~\cite[Corollary 2.4.5]{Springer}),
we know that every element $g\in G$ has a unique decomposition of the form $g= tu$, where $t\in T$ and $u\in U$. 
In this decomposition, if $u=1$ (resp. $t=1$), then $g$ is a semisimple element (resp. a unipotent element).
For $t\in T$, $u\in U$, we define the {\em $t$-twisted $G$-conjugacy class of $u$}, denoted $C_t(u)$, 
by 
\begin{align}\label{A:tconjugacy}
C_t(u) := \{ s ((t^{-1}vt ) u v^{-1})s^{-1} :\ s\in T, v\in U\}.
\end{align}
It is easy to see that $C_t(u)$ is contained in $U$. 
Note also that the $t$-twisted conjugacy classes can be defined more generally for any torus normalizing a unipotent group.

\begin{Theorem}\label{T:pconjugacy}
Let $K$ be a connected solvable linear algebraic group with an abelian unipotent radical $V$.
Let $S$ be a maximal torus of $K$ such that 
$K=S\ltimes V$. 
Let $g$ be an element of $K$ with the Jordan decomposition $g=tu$, where $t\in S$ and $u\in V$.
\begin{enumerate}
\item 
If $g$ is a semisimple element, then there exists a unique element $h\in S$ such that $g$ is conjugate to $h$. 
The $K$-conjugacy class of $h$ is given by its $U$-conjugacy class in $K$, that is, $\{  u h u^{-1}:\  u \in U\}$.
\item If $g$ is a unipotent element, then the $K$-conjugacy class of $g$ is given by its $S$-conjugacy class in $V$.
\item If $g$ is not unipotent or semisimple, then its $K$-conjugacy class is given by $\{ tw' :\ w' \in C_t(u)\}$,
where $C_t(u)$ is the $t$-twisted conjugacy class of $u$ in $V$.
\end{enumerate}
\end{Theorem}

\begin{proof}
The first part of our first claim holds true in a greater generality; in a connected algebraic group any two maximal tori are conjugate to each other~\cite[Theorem 6.4.1]{Springer}. 
In particular, every semisimple element of $K$ is conjugate to a unique element of $S$.
The uniqueness is a consequence of the commutativity of $S$. 
To calculate the $K$-conjugacy class of $h$, let $x$ be an element of $K$.
We write $x$ in the form $x= vs$, where $s\in S$ and $v\in U$. 
Then we have $xhx^{-1} = vs h s^{-1} v^{-1} = vhv^{-1}$. 
This finishes the proof of our first claim.

We proceed to the proof of our second claim.
Let $x\in K$. 
Then $x$ can be written as $x=sv$, where $s\in S$ and $v\in V$. 
Now, if $u$ is from $V$, then we have
\[
xux^{-1} = sv u v^{-1}s^{-1} =sus^{-1} \qquad  \text{since $V$ is abelian}. 
\]
Since $V$ is normalized by $S$, we see that the conjugacy class of $u$ in $K$ is given by the $S$-conjugacy class of $u$ in $V$.  
This finishes the proof of our second assertion.

 Now we end with the proof of the third claim. 
Let $x\in K$. 
such that $x=sv$, where $s\in S$ and $v\in V$, by using the abelian property of the subgroups $S$ and $V$, we have 
\[
xgx^{-1} = sv\ tu\ v^{-1} s^{-1} = t s t^{-1} vt uv^{-1}s^{-1} = t (s w s^{-1}),
\]
where $w= t^{-1}vt u v^{-1} \in V$, hence $sws^{-1}\in C_t(u)$. 
This finishes the proof of our last assertion. 
\end{proof}

We now apply our previous theorem to the $\gh$-classes of $M$. 
Let $X\in M$. 
We will make a practical assumption.
Since every $\gh$-class of $M$ is isomorphic to the unit group of the antichain monoid of a complete bipartite poset with respect to a maximal antichain (Theorem~\ref{T:isomorphicHclasses}), we assume that $H_X$ is given by the unit group of $M$. 
Now let $\mbf{G}$ denote the unit group of $M$. 
Let $\mbf{S}$ denote the maximal diagonal torus of $\mbf{G}$. 
Let $\mbf{V}$ denote the unipotent radical of $\mbf{G}$. 
\medskip

\begin{Corollary}\label{C:moreprecisely}
Let $X\in \mbf{G}$. If $X=TU$ for some $T\in \mbf{S}$ and $U\in \mbf{V}$,
then the $p$-conjugacy class of $X$, denoted $[X]_p$, is given by one of the following items: 
\begin{enumerate}
\item If $U=\mbf{1}_n$, then $[X]_p=\{ U'X U'^{-1} :\ U'\in \mbf{V}\}$.
\item If $T=\mbf{1}_n$, then $[X]_p=\{ T'X T'^{-1} :\ T'\in \mbf{S}\}$.
\item If $U$ and $T$ are different from $\mbf{1}_n$, then $[X]_p=\{TU':\ U'\in C_T(U)\}$.
\end{enumerate}
\end{Corollary}
\begin{proof} 
The proof of Theorem~\ref{T:bipartiteorthodox} shows that $\mbf{V}$ is abelian. 
The rest of the proof follows from a direct application of Theorem~\ref{T:pconjugacy}. We omit the details.  
\end{proof}

The last conjugacy relation we consider here is the $o$-conjugacy relation which is defined by Otto in~\cite{Otto}. 
Let $a$ and $b$ be two elements from a semigroup $S$. 
Then we say that $a$ and $b$ are {\em $o$-conjugacy related}, and write $a\simo b$ if there exists $\{g,h\}\subseteq S^1$
such that 
\begin{align*}
ag = gb \qquad\text{and}\qquad bh = ha.
\end{align*}
In fact, it is easy to verify that $\simo$ is an equivalence relation. 
It is known that $\simp\ \subseteq\ \simo$ and that if $S$ has a zero, then $\simo$ is a universal relation. 
Since antichain monoids do not have a zero element, it is interesting to know the properties of $\simo$ on antichain monoids.
It turns out that $\simo$ is still a universal relation on the idempotent semigroup of $M$. 

\begin{Proposition}\label{P:allidempotentsoconjugate}
Let $M$ be an antichain monoid as before. 
If $X$ and $Y$ are two elements of $E(M)$, then we have $X\simo Y$ in $E(M)$.
\end{Proposition}

\begin{proof}
Let $X$ and $Y$ be given by the matrices $\begin{bmatrix} \mbf{1}_k & B_X \\ 0 &  \diag_X \end{bmatrix}$ and 
$\begin{bmatrix} \mbf{1}_k & B_Y \\ 0 &  \diag_Y \end{bmatrix}$, respectively. 
Let $Z$ denote the element of $M$ defined by 
$\begin{bmatrix} \mbf{1}_k & B_Y \\ 0 &  0 \end{bmatrix}$.
Clearly, $Z$ is an idempotent. 
Then we have 
\[
XZ = 
\begin{bmatrix} \mbf{1}_k & B_Y \\ 0 &  0\end{bmatrix}\qquad\text{and}\qquad
ZY =
\begin{bmatrix} \mbf{1}_k & B_Y+B_Y\diag_Y \\ 0 &  0\end{bmatrix}.
\]
But the identity $Y^2 = Y$ implies that $B_Y+B_Y \diag_Y = B_Y$, hence that $B_Y\diag_Y=0$. 
In other words, we have $XZ=ZY$. 
A similar argument shows that if we choose $W\in E(M)$ as 
$W=\begin{bmatrix} \mbf{1}_k & B_X \\ 0 &  0 \end{bmatrix}$, then we have $YW = WX$.
This finishes the proof of our proposition.
\end{proof}

\begin{Problem}
Theorem~\ref{T:modifiedlastannounced} solves the problem of $p$-conjugacy only in the case of antichain monoids of (complete) bipartite posets. The general case remains open. 
Also unsolved is the equivalence of  $\simp$, $\simn$ and $\simtr$  on an arbitrary incidence monoid. 
It would be interesting to investigate the structures of the individual conjugacy classes under these relations.

A closely related question is about the determination of the $\simo$-conjugacy classes in $M$. 
Although Proposition~\ref{P:allidempotentsoconjugate} shows that $\simo$ is a universal relation within $E(M)$,
in general, $M$ has many distinct $\simo$-conjugacy classes. 
Is it true that $E(M)$ is a single $\simo$-conjugacy class in $M$? 
\end{Problem}

\subsection*{Acknowledgements} 
The authors thank the Centre for Mathematics and Applications of the NOVA University of Lisbon for making this collaboration possible. 
The first author is partially supported by a grant from the Louisiana Board of Regents (090ENH-21).  The second and third authors were supported by national funds through the FCT -- Funda\c{c}\~{a}o para a Ci\^{e}ncia e a Tecnologia, I.P., under the scope of the projects {\scshape UIDB}/00297/2020 and {\scshape UIDP}/00297/2020 (Center for Mathematics and Applications) and {\scshape PTDC}/{\scshape MAT-PUR}/31174/2017.
The authors thank Michel Brion for his valuable comments. 
Finally, the authors thank the referee whose comments and suggestions improved the quality of the paper.

\subsection*{Ethics declarations/Conflict of interests}
The authors have no conflicts of interest to declare that are relevant to the content of this article.

\subsection*{Data availability statement}
Data sharing not applicable to this article as no datasets were generated or analyzed during the current study.

\bibliography{references.bib}
\bibliographystyle{plain}
\end{document}